%% file: main.tex
\documentclass[journal,twoside,web]{ieeecolor}
\usepackage{generic}
\usepackage{textcomp}

\usepackage{graphicx}
\usepackage{times}
\usepackage{amsmath,amsbsy,amssymb,mathtools}
\usepackage{color}
\usepackage[us]{datetime}
\usepackage{url}
\usepackage{theorem}
\usepackage{multirow}
\usepackage{footnote}
\usepackage{nomencl}
\usepackage{cite}
\usepackage{algorithm,tabularx}
\usepackage[noend]{algpseudocode}
\usepackage{stackengine}

\usepackage{commath}
\usepackage{mathrsfs}
\usepackage{setspace}

\newcommand{\negpar}[1][-1em]{%
  \ifvmode\else\par\fi
  {\parindent=#1\leavevmode}\ignorespaces
}
\usepackage{stackengine}
\setstackEOL{\\}
\makeatletter
\def\BState{\State\hskip-\ALG@thistlm}
\makeatother

\newcommand{\R}{\mathbb{R}} % Real number
\theoremstyle{plain}
\newtheorem{theorem}{Theorem}

\newtheorem{lemma}{Lemma}
\newtheorem{corollary}{Corollary}
\newtheorem{remark}{Remark}
\newtheorem{problem}{Problem}
\newtheorem{assumption}{Assumption}

\usepackage{lipsum}
\newcommand{\qed}{\hfill\rule{1.5ex}{1.5ex}}

\algnewcommand\algorithmicswitch{\textbf{switch}}
\algnewcommand\algorithmiccase{\textbf{case}}
\algdef{SE}[SWITCH]{Switch}{EndSwitch}[1]{\algorithmicswitch\ #1\ \algorithmicdo}{\algorithmicend\ \algorithmicswitch}%
\algdef{SE}[CASE]{Case}{EndCase}[1]{\algorithmiccase\ #1}{\algorithmicend\ \algorithmiccase}%
\algtext*{EndSwitch}%
\algtext*{EndCase}%

\def\BibTeX{{\rm B\kern-.05em{\sc i\kern-.025em b}\kern-.08em
    T\kern-.1667em\lower.7ex\hbox{E}\kern-.125emX}}

\begin{document}

\onecolumn
\copyright 2021 IEEE. Personal use of this material is permitted. Permission from IEEE must be obtained for all other uses, in any current or future media, including reprinting/republishing this material for advertising or promotional purposes, creating new collective works, for resale or redistribution to servers or lists, or reuse of any copyrighted component of this work in other works.

This work has been submitted to the IEEE for possible publication. Copyright may be transferred without notice, after which this version may no longer be accessible.
\twocolumn
\newpage

\title{Hamilton-Jacobi Equations for Two Classes of State-Constrained 
Zero-Sum Games}
% \title{\large \bf
% % \title{
% Hamilton-Jacobi Equations for Two Classes of State-Constrained 
% Zero-Sum Games
% % Optimal Control Problems and Zero-Sum Games %\\Part I. Hamilton-Jacobi Equations
% }
\author{Donggun Lee, \IEEEmembership{Student Member, IEEE}, and Claire J. Tomlin, \IEEEmembership{Fellow, IEEE} % <-this % stops a space
\thanks{This research is supported by ONR under the BRC program in multibody control systems, by DARPA under the Assured Autonomy program, and by
NSF grant \#1837244.}% <-this % stops a space
\thanks{Donggun Lee is with the Department of Mechanical Engineering, University of California, Berkeley, USA.
        {\tt\small donggun\_lee@berkeley.edu}}%
\thanks{Claire J. Tomlin is with the Department of Electrical Engineering and Computer Sciences, University of California, Berkeley, USA.
        {\tt\small tomlin@eecs.berkeley.edu}}%
% \thanks{$^{*}$ Both authors contributed equally to this work.}%
}

\maketitle
% \thispagestyle{empty}
% \pagestyle{empty}
% \thispagestyle{plain}
% \pagestyle{plain}

%%%%%%%%%%%%%%%%%%%%%%%%%%%%%%%%%%%%%%%%%%%%%%%%%%%%%%%%%%%%%%%%%%%%%%%%%%%%%%%%
\begin{abstract}
This paper presents Hamilton-Jacobi (HJ) formulations for two classes of two-player zero-sum games: one with a maximum cost value over time, and one with a minimum cost value over time.
In the zero-sum game setting, player A minimizes the given cost while satisfying state constraints, and player B wants to prevent player A's success.
For each class of problems, this paper presents two HJ equations: one for time-varying dynamics, cost, and state constraint; the other for time-invariant dynamics, cost, and state constraint.
Utilizing the HJ equations, the optimal control for each player is analyzed, and a numerical algorithm is presented to compute the solution to the HJ equations.
A two-dimensional water system is introduced as an example to demonstrate the proposed HJ framework.
\end{abstract}

%%%%%%%%%%%%%%%%%%%%%%%%%%%%%%%%%%%%%%%%%%%%%%%%%%%%%%%%%%%%%%%%%%%%%%%%%%%%%%%%
\section{Introduction}

In two-player zero-sum games, one player's control signal minimizes a cost while satisfying a state constraint, while the second player's control signal tries either to maximize the cost or to violate the state constraint.
An optimal control problem may be considered a special case of the zero-sum game: a control signal that minimizes the given cost while satisfying the constraint is to be determined.
% , in which the second player does not contribute the system at all.
% In other words, the optimal control problem is to determine 

The Hamilton-Jacobi (HJ) partial differential equation (PDE) can be used to represent zero-sum games for dynamical systems.
The HJ formulation includes a cost function in the form of the integration of a stage cost and a terminal cost, nonlinear dynamics, control constraints, and state constraints.

The zero-sum game can be classified according to 1) whether the terminal time is a given constant or a variable to be determined, and 2) whether or not state constraints exist.
If the problem is state-unconstrained and the terminal time is a given constant, the Hamilton-Jacobi-Isaacs (HJI) PDE \cite{evans1984differential} applies.
For the state-constrained problem where the terminal time is given, \cite{altarovici2013general} presents the corresponding HJ equation. 
For problems where the terminal time is a variable to be determined, \cite{mitchell2005time} deals with the state-unconstrained problem where the stage cost is zero, and \cite{margellos2011hamilton,fisac2015reach} deal with the zero-stage-cost and state-constrained problems.
This paper generalizes the previous work to deal with the case of non-zero stage-cost and state constraint.

This paper proposes HJ equations for two classes of state-constrained problems where the terminal time is a variable to be determined and the stage cost is non-zero.
In the first class of problems, player A wants to minimize the maximum cost over time while satisfying the state constraint, and player B wants to prevent player A's success.
This class of problems can be interpreted as a robust control problem on the time and disturbances that optimizes the maximum cost over time with respect to the worst disturbances.
In the second class of problems, player A wants to minimize the minimum cost over time while satisfying the state constraint, and player B again wants to prevent player A's success.
This class of problems can be interpreted as another robust control problem on the disturbances that optimizes the minimum cost over time with respect to the worst disturbances.

The proposed HJ equations can generally deal with both time-varying and time-invariant dynamics, cost, and state constraint.
Furthermore, this paper presents additional HJ equations equivalent to the proposed HJ equations for the time-invariant case.

\subsection{Contribution}

This paper presents four HJ equations for zero-sum games: two classes, and time-varying and time-invariant.
Among the four HJ equations, three equations are proposed by this paper, and the other one is presented in \cite{donggun2020CDC} and reviewed here for completeness.
Also, this paper provides and presents analysis for the optimal control signal for each player, numerical algorithm to compute the proposed HJ equation, as well as a practical example.

\subsection{Organization}

The organization of this paper is as follows. 
Section \ref{sec:problemDefinition} presents a mathematical formulation for two classes of state-constrained zero-sum games.
Section \ref{sec:HJeq_game} presents the HJ equations for the first class of problems both time-varying and time-invariant.
Section \ref{sec:HJeq_Prob2Big} presents the HJ equations for the second class of problems for both the time-varying and time-invariant cases.
Section \ref{sec:optCtrl_strategy} presents analysis for an optimal control signal based on the solution to the HJ equations.
Section \ref{sec:NumericalAlgorithm} presents a numerical algorithm to compute the solution to the HJ equations for each class of problem.
Section \ref{sec:example} provides a practical example where our HJ formulation can be utilized, and Section \ref{sec:conclusion} concludes this paper.
Proofs are detailed in the Appendices.

%%%%%%%%%%%%%%%%%%%%%%%%%%%%%%%%%%%%%%%%%%%%%%%%%%%%%%%%%%%%%%%%%%%%%%%%%%%%%%%%%%%%%%%%%%%%%%%%%%%%%%%%%%%%%%%%%%%%%%%%%%%%%%%%%%%%%%%%%%%%%%%%%%%%%%%%%%%%%%%%%%%%%%%%%%%%%%%%%%%%%%%%%%%%%%%%%%%%%%%%%%%%%%%%%%%%%%%%%%%%%%%%%%%%%%%%%%%%%%%%%%%%%%%%%%%%%%%%%%%%%%%%%%%%%%%%%%%%%%%%%%%%%%%%%%%%%%%%%%%%%%%%%%%%%%%%
\section{Two Classes of State-Constrained Zero-Sum Games }
\label{sec:problemDefinition}

We first present the two classes of two-player zero-sum games, called Problems \ref{prob:1} and \ref{prob:2}.
Consider a dynamical system:
\begin{align}
    \dot{\mathrm{x}}(s) = f(s,\mathrm{x}(s),\alpha(s),\beta(s)), s\in[t,T], \text{  and  } \mathrm{x}(t)=x,
    \label{eq:dynamics}
\end{align}
where $(t,x)$ are the initial time and state, $\mathrm{x}:[t,T]\rightarrow\R^n$ is the state trajectory, $f:[0,T]\times \R^n\times A\times B\rightarrow\R^n$ is the dynamics, $A\subset \R^{m_a},B\subset \R^{m_b}$ are the control constraints,  $\alpha\in\mathcal{A}(t),\beta\in\mathcal{B}(t)$ are the control signals, in each, player A controls $\alpha$ and player B controls $\beta$, and the sets of measurable control signals are
\begin{align}
\begin{split}
    &\mathcal{A}(t) \coloneqq \{\alpha:[t,T]\rightarrow A ~|~ \|\alpha\|_{L^\infty(t,T)} <\infty\},\\
    &\mathcal{B}(t) \coloneqq \{\beta:[t,T]\rightarrow B ~|~ \|\beta\|_{L^\infty(t,T)} <\infty\}.
\end{split}
    \label{eq:def_CtrlTraj}
\end{align}

In each zero-sum game, we specify each player's control signal or strategy: player A wants to minimize the cost under the state constraint, and player B wants to prevent player A's success, although the cost is defined in different ways for Problems \ref{prob:1} and \ref{prob:2}.
In each problem, we introduce two value functions depending on which players play first or second.

\begin{problem} For given initial time and state $(t,x)$, solve 
\begin{align}
    \begin{split}
    & \vartheta_1^+ (t,x) \coloneqq \sup_{\delta\in\Delta(t)}\inf_{\alpha\in\mathcal{A}(t)} \\
    &\quad \max_{\tau\in[t,T]} \int_t^\tau L(s,\mathrm{x}(s),\alpha(s),\delta[\alpha](s))ds + g(\tau,\mathrm{x}(\tau)), 
    \label{eq:def_vartheta1_upper}
    \end{split}
    \\
    & \text{subject to } \quad c(s,\mathrm{x}(s)) \leq 0, \quad s\in [t,\tau],
    \label{eq:def_vartheta_const1_upper}
\end{align}
where $\mathrm{x}$ solves \eqref{eq:dynamics} for $(\alpha,\delta[\alpha])$; and solve 
\begin{align}
    \begin{split}
    & \vartheta_1^- (t,x) \coloneqq \inf_{\gamma\in\Gamma(t)}\sup_{\beta\in\mathcal{B}(t)} \\
    &\quad \max_{\tau\in[t,T]} \int_t^\tau L(s,\mathrm{x}(s),\gamma[\beta](s),\beta(s))ds + g(\tau,\mathrm{x}(\tau)), 
    \label{eq:def_vartheta1_lower}
    \end{split}
    \\
    & \text{subject to } \quad c(s,\mathrm{x}(s)) \leq 0, \quad s\in [t,\tau],
    \label{eq:def_vartheta_const1_lower}
\end{align}
where $\mathrm{x}$ solves \eqref{eq:dynamics} for $(\gamma[\beta],\beta)$.
\label{prob:1}
\end{problem}
$\Delta(t)$ is a set of non-anticipative strategies for player B, and $\Gamma(t)$ is a set of non-anticipative strategies for player A.
The non-anticipative strategy outputs a control signal for the second player as a reaction to the first player's control signal without using the future information.
The non-anticipative strategy has been introduced by Elliott and Kalton \cite{elliott1972existence}:
% , which is a reaction to player A's control signal.
\begin{align}
    \Delta(t)\coloneqq&\{\delta:\mathcal{A}(t)\rightarrow\mathcal{B}(t) ~| ~ \forall s\in[t,\tau] \text{ and } \alpha,\bar{\alpha}\in\mathcal{A}(t),\notag\\
    &\text{      if } \alpha(\tau)=\bar{\alpha}(\tau) \text{ a.e. } \tau\in[t,s],\notag\\
    &\text{      then } \delta[\alpha](\tau)=\delta[\bar{\alpha}](\tau) \text{ a.e. } \tau\in[t,s] \}.
\label{eq:playerB_strategySet}\\
\Gamma(t)\coloneqq&\{\gamma:\mathcal{B}(t)\rightarrow\mathcal{A}(t) ~| ~ \forall s\in[t,\tau], \beta,\bar{\beta}\in\mathcal{B}(t),\notag\\
    &\text{      if } \beta(\tau)=\bar{\beta}(\tau) \text{ a.e. } \tau\in[t,s],\notag\\
    &\text{      then } \gamma[\beta](\tau)=\gamma[\bar{\beta}](\tau) \text{ a.e. } \tau\in[t,s] \}.
% \end{split}
\label{eq:playerA_strategySet}
\end{align}

\begin{problem} For given initial time and state $(t,x)$, solve 
\begin{align}
    \begin{split}
    & \vartheta_2^+ (t,x) \coloneqq \sup_{\delta\in\Delta(t)}\inf_{\alpha\in\mathcal{A}(t)} \\
    &\quad \min_{\tau\in[t,T]} \int_t^\tau L(s,\mathrm{x}(s),\alpha(s),\delta[\alpha](s))ds + g(\tau,\mathrm{x}(\tau)), 
    \label{eq:def_vartheta2_upper}
    \end{split}
    \\
    & \text{subject to } \quad c(s,\mathrm{x}(s)) \leq 0, \quad s\in [t,\tau],
    \label{eq:def_vartheta_const2_upper}
\end{align}
where $\mathrm{x}$ solves \eqref{eq:dynamics} for $(\alpha,\delta[\alpha])$; and solve 
\begin{align}
    \begin{split}
    & \vartheta_2^- (t,x) \coloneqq \inf_{\gamma\in\Gamma(t)}\sup_{\beta\in\mathcal{B}(t)} \\
    &\quad \min_{\tau\in[t,T]} \int_t^\tau L(s,\mathrm{x}(s),\gamma[\beta](s),\beta(s))ds + g(\tau,\mathrm{x}(\tau)), 
    \label{eq:def_vartheta2_lower}
    \end{split}
    \\
    & \text{subject to } \quad c(s,\mathrm{x}(s)) \leq 0, \quad s\in [t,\tau],
    \label{eq:def_vartheta_const2_lower}
\end{align}
where $\mathrm{x}$ solves \eqref{eq:dynamics} for $(\gamma[\beta],\beta)$.
\label{prob:2}
\end{problem}
For both problems, $L:[t,T]\times\R^n\times A\times B\rightarrow \R$ is the stage cost, $g:\R\times\R^n\rightarrow\R$ is the terminal cost, $f:[t,T]\times\R^n\times A \times B\rightarrow\R^n$ is the system dynamics, and $c:[t,T]\times\R^n\rightarrow\R$ is the state constraint function.

The difference between $\vartheta_i^+$ and $\vartheta_i^-$ ($i=1,2$) is play order.
In $\vartheta_i^+(t,x)$, at each time $s\in[t,T]$, player A first plays $\alpha(s)$, and then player B reacts by following its own strategy $\delta[\alpha](s)$.
Despite this play order at each time, the choice of player B's strategy comes first since it should be chosen without information about player A's control signal.
In other words, player B first chooses its strategy, and then player A chooses its control signal.
In $\vartheta_i^-(t,x)$, at each time $s$, player B first plays $\beta(s)$, and then player A reacts with its strategy $\gamma[\beta](s)$.
Similarly to $\vartheta_i^+(t,x)$, in $\vartheta_i^-(t,x)$, player A first chooses its strategy, and then player B chooses its control signal.

Problems \ref{prob:1} and \ref{prob:2} are representative of many practical problems.
For Problem \ref{prob:1}, consider two water systems where player A controls the water level of pond 1 that is connected to pond 2.
Suppose player B is precipitation. 
Player A needs to minimize the highest water level of pond 1 over time while satisfying constraints for water level of pond 1 and 2 under the worst precipitation assumption.
For Problem \ref{prob:2}, consider a car that tries to change its lane while avoiding collision with other cars. Here, the cost is the distance to the goal lane, and the car wants to successfully change lanes at some time in the given time interval, while other cars might bother to lane change.

This paper assumes the following.
\begin{assumption}[Lipschitz continuity and compactness]
~
\begin{enumerate}
    \item $A$ and $B$ are compact;
    
    \item $f:[0,T]\times\R^n\times A\times B\rightarrow\R^n$, $f=f(t,x,a,b)$ is Lipschitz continuous in $(t,x)$ for each $(a,b)\in A\times B$;

    \item the stage cost $L:[0,T]\times\R^n\times A\times B \rightarrow\R$, $L=L(t,x,a,b)$ is Lipschitz continuous in $(t,x)$ for each $(a,b)\in A\times B$;

    \item for all $(t,x)\in[0,T]\times\R^n$, $\{f(t,x,a,b)~|~a\in A, b\in B\}$ and $\{L(t,x,a,b)~|~a\in A, b\in B\}$ are compact and convex;
    
    \item the terminal cost $g:[0,T]\times\R^n\rightarrow\R$, $g=g(t,x)$ is Lipschitz continuous in $(t,x)$;

    \item the state constraint $c:[0,T]\times\R^n \rightarrow\R$, $c=c(t,x)$ is Lipschitz continuous in $(t,x)$;

    \item the stage cost ($L$) and the terminal cost ($g$) are bounded below.
    
\end{enumerate}
\label{assum:BigAssum}
\end{assumption}

\section{Hamilton-Jacobi Equations for Problem \ref{prob:1}}
\label{sec:HJeq_game}

\subsection{HJ equation for Problem 1 (time-varying case)}
\label{sec:HJeq_gameProb1}

In this subsection, we derive an HJ equation for Problem \ref{prob:1} ($\vartheta_1^\pm$).
Unfortunately, for some initial time and state $(t,x)$, there is no control $\alpha$ (or strategy $\gamma$) of player A that satisfies the state constraint for all strategies $\delta$ of player B (or control signal $\beta$). In this case, $\vartheta_1^\pm(t,x)$ is infinity. Thus, $\vartheta_1^\pm$ is neither continuous nor differentiable in $(0,T)\times\R^n$.
% , in which it is hard to deal with differential equations.

To overcome this issue, we utilize an additional variable $z\in\R$ to define continuous value functions $V_1^\pm$ in \eqref{eq:game_V1_upper} and \eqref{eq:game_V1_lower} that combine the cost $\vartheta_1^\pm$ in \eqref{eq:def_vartheta1_upper} or \eqref{eq:def_vartheta1_lower}, and the constraint in \eqref{eq:def_vartheta_const1_upper} or \eqref{eq:def_vartheta_const1_lower}. 
We call this method \textit{the augmented-$z$ method}. 
This method has been utilized to handle state constraints to solve other HJ problems \cite{altarovici2013general,donggun2020CDC}.
$V_1^\pm$ is well-defined in $[0,T]\times\R^n\times\R$. 
\begin{align}
    V_1^+(t,x,z) \coloneqq    \sup_{\delta\in\Delta(t)}\inf_{\alpha\in\mathcal{A}(t)} J_1(t,x,z,\alpha,\delta[\alpha]),
    \label{eq:game_V1_upper}\\
% \end{align}
% and 
% \begin{align}
    V_1^-(t,x,z) \coloneqq    \inf_{\gamma\in\Gamma(t)}\sup_{\beta\in\mathcal{B}(t)} J_1(t,x,z,\gamma[\beta],\beta),
    \label{eq:game_V1_lower}
\end{align}
where cost $J_1:(t,x,z,\alpha,\beta)\rightarrow\R$ is defined as follows:
\begin{align}
\begin{split}
    &J_1(t,x,  z,\alpha,\beta) \coloneqq \max_{\tau\in[t,T]}\max\bigg\{
    \max_{s\in[t,\tau]}c(s,\mathrm{x}(s)), \\
    &\quad\quad \int_t^\tau L(s,\mathrm{x}(s),\alpha(s),\beta(s))ds+g(\tau,\mathrm{x}(\tau))-z
    \bigg\},
    \end{split}
    \label{eq:ctrl_costDef1}
\end{align}
where $\mathrm{x}$ solves \eqref{eq:dynamics}.
% For the auxiliary variable $z$, 
Define the auxiliary state trajectory $\mathrm{z}$ solving
\begin{align}
    \dot{\mathrm{z}}(s) = -L(s,\mathrm{x}(s),\alpha(s),\beta(s)), s\in[t,T], \text{ and } \mathrm{z}(t)=z.
    \label{eq:dynamics_z}
\end{align}
Then, \eqref{eq:dynamics} and \eqref{eq:dynamics_z} are the joint ODEs whose solution is the augmented state trajectories: $(\mathrm{x},\mathrm{z}):[t,T]\rightarrow\R^{n+1}$
\begin{align}
\begin{small}
    \begin{bmatrix}\dot{\mathrm{x}}(s)\\\dot{\mathrm{z}}(s)\end{bmatrix} = \begin{bmatrix}f(s,\mathrm{x}(s),\alpha(s),\beta(s))\\-L(s,\mathrm{x}(s),\alpha(s),\beta(s))\end{bmatrix}, s\in[t,T], \begin{bmatrix}\mathrm{x}(t)\\\mathrm{z}(t)\end{bmatrix} =\begin{bmatrix}x\\z\end{bmatrix}.
    \end{small}
    \label{eq:dynamics_comb}
\end{align}
Then, $J_1$ in \eqref{eq:ctrl_costDef1} becomes
\begin{align}
\begin{split}
    J_1=& \max_{\tau\in[t,T]}\max\Big\{
    \max_{s\in[t,\tau]}c(s,\mathrm{x}(s)),  g(\tau,\mathrm{x}(\tau))-\mathrm{z}(\tau)
    \Big\} \\
    =& \max\Big\{
    \max_{s\in[t,T]}c(s,\mathrm{x}(s)), \max_{\tau\in[t,T]} g(\tau,\mathrm{x}(\tau))-\mathrm{z}(\tau)
    \Big\}.
    \end{split}
    \label{eq:ctrl_costDef1_2}
\end{align}
The last equality is derived by the distributive property of the maximum operations.

Lemma \ref{lemma:Equiv_twoCost1} shows that $\vartheta_1^\pm$ can be found if $V_1^\pm$ are known. 
For initial time and state $(t,x)$ for which there is no control or strategy of player A such that the state constraint ($c(s,\mathrm{x}(s))\leq 0,s\in[t,T]$) is satisfied for player B's best control signal or strategy, $V_1^\pm(t,x,z)$ is always greater than 0 for all $z\in\R$. In this case, Lemma \ref{lemma:Equiv_twoCost1} implies that $\vartheta_1^\pm(t,x)$ is infinity.

\begin{lemma}[Equivalence of two value functions] Suppose Assumption \ref{assum:BigAssum} holds. For all $(t,x)\in[0,T]\times\R^n$, $\vartheta_1^+$ (\eqref{eq:def_vartheta1_upper} subject to \eqref{eq:def_vartheta_const1_upper}), $\vartheta_1^-$ (\eqref{eq:def_vartheta1_lower} subject to \eqref{eq:def_vartheta_const1_lower}), $V_1^+$ in \eqref{eq:game_V1_upper}, and $V_1^-$ in \eqref{eq:game_V1_lower} have the following relationship.
    \begin{align}
    % \begin{split}
        \vartheta_1^\pm(t,x) = \min z \text{ subject to } V_1^\pm(t,x,z)\leq 0.
        % \end{split}
        \label{eq:lemma_equiv_twoCost1}
    \end{align}
This implies that 
    \begin{align}
    \begin{split}
     \vartheta_1^+ (t,x) &= \sup_{\delta\in\Delta(t)}\inf_{\alpha\in\mathcal{A}(t)} \max_{\tau\in[t,T]} \\
    & \int_t^\tau L(s,\mathrm{x}(s),\alpha(s),\delta[\alpha](s))ds + g(\tau,\mathrm{x}(\tau)), 
    \label{eq:vartheta1_upper_eq}
    \end{split}
    \\
    & \text{subject to }\quad  c(s,\mathrm{x}(s)) \leq 0, \quad s\in [t,T],
    \label{eq:vartheta1_upper_eq_const}
    \end{align}
    where $\mathrm{x}$ solves \eqref{eq:dynamics} for $(\alpha,\delta[\alpha])$, and 
    \begin{align}
    \begin{split}
     \vartheta_1^- (t,x) &= \inf_{\gamma\in\Gamma(t)}\sup_{\beta\in\mathcal{B}(t)} \max_{\tau\in[t,T]} \\
    & \int_t^\tau L(s,\mathrm{x}(s),\gamma[\beta](s),\beta(s))ds + g(\tau,\mathrm{x}(\tau)), 
    \label{eq:vartheta1_lower_eq}
    \end{split}
    \\
    & \text{subject to }\quad  c(s,\mathrm{x}(s)) \leq 0, \quad s\in [t,T],
    \label{eq:vartheta1_lower_eq_const}
    \end{align}
    where $\mathrm{x}$ solves \eqref{eq:dynamics} for $(\gamma[\beta],\beta)$.
\label{lemma:Equiv_twoCost1}
\end{lemma}
\textbf{Proof.} See Appendix \ref{appen:lemma_Equiv_twoCost1}.

The rest of this subsection focuses on the derivation of the corresponding HJ equation for $V_1^\pm$.
The HJ equation is based on the principle of dynamic programming in Lemma \ref{lemma:dynamicProgramming}.
\begin{lemma}[Optimality condition]
    Fix $(t,x,z)\in[0,T]\times\R^n\times\R$.  Consider a small step $h>0$ such that $t+h\leq T$, $V_1^+$ \eqref{eq:game_V1_upper} has the following property:
    \begin{align}
    \begin{split}
        &V_1^+ (t,x,z) = \sup_{\delta\in\Delta(t)}\inf_{\alpha\in\mathcal{A}(t)}\max\Big\{\max_{s\in[t,t+h]} c(s,\mathrm{x}(s)) , \\ 
        &  \max_{s\in[t,t+h]} g(\mathrm{x}(s))-\mathrm{z}(s),
         V_1^+ (t+h,\mathrm{x}(t+h),\mathrm{z}(t+h)) \Big\},
    \end{split}
    \label{eq:dynamicProgramming_V1_upper}
    \end{align}
    where $(\mathrm{x},\mathrm{z})$ solves \eqref{eq:dynamics_comb} for $(\alpha,\delta[\alpha])$.
    Similarly, for $V_1^-$ \eqref{eq:game_V1_lower},
    \begin{align}
    \begin{split}
        &V_1^- (t,x,z) = \inf_{\gamma\in\Gamma(t)}\sup_{\beta\in\mathcal{B}(t)}\max\Big\{\max_{s\in[t,t+h]} c(s,\mathrm{x}(s)) , \\ 
        &  \max_{s\in[t,t+h]} g(\mathrm{x}(s))-\mathrm{z}(s),
         V_1^- (t+h,\mathrm{x}(t+h),\mathrm{z}(t+h)) \Big\},
    \end{split}
    \label{eq:dynamicProgramming_V1_lower}
    \end{align}
    where $(\mathrm{x},\mathrm{z})$ solves \eqref{eq:dynamics_comb} for $(\gamma[\beta],\beta)$.
    \label{lemma:dynamicProgramming}
\end{lemma}
\textbf{Proof.} See Appendix \ref{appen:lemma_dynamicProgramming}.

Theorem \ref{thm:HJeq_prob1} presents the corresponding HJ equations for $V_1^\pm$ in \eqref{eq:game_V1_upper} and \eqref{eq:game_V1_lower} using viscosity theory.
Intuitively, the HJ equation in Theorem \ref{thm:HJeq_prob1} is derived as $h$ in Lemma \ref{lemma:dynamicProgramming} converges to zero.

\begin{theorem}\textnormal{\textbf{(HJ equation for Problem \ref{prob:1})}}
    For all $(t,x,z)\in[0,T]\times\R^n\times\R$, $V_1^\pm$ in \eqref{eq:game_V1_upper} and \eqref{eq:game_V1_lower} is the unique viscosity solution to the HJ equation:
    \begin{align}
    \begin{split}
        \max\Big\{c(t,x)-&V_1^\pm(t,x,z), g(t,x)-z-V_1^\pm(t,x,z), \\ 
        &V_{1,t}^\pm - \bar{H}^\pm (t,x,z,D_x V_1^\pm,D_z V_1^\pm)\Big\}=0
        \end{split}
        \label{eq:HJeq1}
    \end{align}
    in $(0,T)\times\R^n\times\R$, where $\bar{H}^\pm:[0,T]\times\R^n\times\R\times\R^n\times\R\rightarrow\R$ \begin{align}
        \bar{H}^+(t,x,z,p,q) \coloneqq \max_{a\in A}\min_{b\in B} - p\cdot f(t,x,a,b) + q L(t,x,a,b),
        \label{eq:def_Hambar_+}\\
        \bar{H}^-(t,x,z,p,q) \coloneqq \min_{b\in B}\max_{a\in A} - p\cdot f(t,x,a,b) + q L(t,x,a,b),
        \label{eq:def_Hambar_-}
    \end{align}
    and
    \begin{align}
        V_1^\pm(T,x,z) = \max\{ c(T,x), g(T,x)-z \}
        \label{eq:V_terminalValue1}
    \end{align}
    on $\{t=T\}\times\R^n\times\R$. Denote $V_{1,t}^\pm=\frac{\partial V_1^\pm}{\partial t}$, $D_x V_1^\pm = \frac{\partial V_1^\pm}{\partial x}$, and $D_z V_1^\pm = \frac{\partial V_1^\pm}{\partial z}$.
    \label{thm:HJeq_prob1}
\end{theorem}
\textbf{Proof.} See Appendix \ref{appen:thm_HJeq_prob1}.

\subsection{HJ equation for Problem 1 (time-invariant case)}
\label{sec:timeInvProb}

We define the problem as time-invariant if the stage cost, terminal cost, dynamics, and state constraints are all independent of time.

In this section, we convert $\vartheta_1^\pm$ (\eqref{eq:def_vartheta1_upper} subject to \eqref{eq:def_vartheta_const1_upper} and \eqref{eq:def_vartheta1_lower} subject to \eqref{eq:def_vartheta_const1_lower}) to a fixed-terminal-time problem for the time-invariant case of Problem \ref{prob:1}, which allows to utilize methods for the fixed-terminal-time problems \cite{altarovici2013general}.
In the fixed-terminal-time problem, optimal control signals of players have to be determined, but the terminal time does not need to be specified but is given.

The conversion of Problem \ref{prob:1} to a fixed-terminal-time problem by introducing a freezing control signal $\mu:[t,T]\rightarrow[0,1]$ to the dynamics and a set of freezing control signals: 
\begin{align}
    \dot{\mathrm{x}}(s) = f(\mathrm{x}(s),\alpha(s),\beta(s))\mu(s), s\in[t,T], \mathrm{x}(t)=x,
    \label{eq:dynamics_timeInvariant}\\
% \end{align}
% Define a set of freezing control signals
% \begin{align}
    \mathcal{M} (t) \coloneqq \{ \mu: [t,T]\rightarrow [0,1]~|~
    \|\mu\|_{L^\infty (t,T)}<\infty\}.
    \label{eq:ctrlSet_freezing}
\end{align}
This freezing control signal controls the contribution of the two players to the system.
For example, $\mu(s)=0$ implies that the state stops at $s$, and the two players do not contribute to the system.
On the other hand, $\mu(s)=1$ allows the state evolves by the control signals of the players. 
% creates another time line that flows faster than the actual time line or skips it at each time.
The maximum over $\tau$ operation in Problem \ref{prob:1} can be replaced by the maximum over the freezing control signal if it eliminates contribution of the two players after the maximal terminal time.

% Using these definitions, w
We present fixed-terminal-time problems as below:
\begin{align}
 \begin{split}
%\sup_{\substack{\delta\in \Delta(t) \\ \nu_A\in \Delta_d(t)}}
    \tilde\vartheta_1^+ (t&,x) \coloneqq \sup_{\delta\in \Delta(t),\nu_A\in \textrm{N}_A(t)} \inf_{\alpha\in\mathcal{A}(t)} \\
    &\int_t^T L(\mathrm{x}(s),\alpha(s),\delta[\alpha](s)) \nu_A[\alpha](s)ds + g(\mathrm{x}(T)),
    \label{eq:vartheta1_upper_opt_cost}
\end{split}\\
&\quad\quad\quad \text{subject to } c(\mathrm{x}(s))\leq 0, s\in[t,T],
\label{eq:vartheta1_upper_opt_const}         
\end{align}
where $\mathrm{x}$ solves \eqref{eq:dynamics_timeInvariant} for $(\alpha,\delta[\alpha],\nu_A[\alpha])$;
\begin{align}
\begin{split}
    \tilde\vartheta_1^- (t&,x) \coloneqq\inf_{\tilde\gamma\in \tilde{\Gamma}(t)} \sup_{\beta\in\mathcal{B}(t),\mu\in\mathcal{M}(t)} \\
    &\int_t^T L(\mathrm{x}(s),\tilde\gamma[\beta,\mu](s),\beta(s)) \mu(s)ds + g(\mathrm{x}(T)),
    \label{eq:vartheta1_lower_opt_cost}
\end{split}\\
&\quad\quad\quad \text{subject to } c(\mathrm{x}(s))\leq 0, s\in[t,T],
\label{eq:vartheta1_lower_opt_const}
\end{align}
where $\mathrm{x}$ solves \eqref{eq:dynamics_timeInvariant} for $(\tilde\gamma[\beta,\mu],\beta,\mu)$.
$\textrm{N}_A$ is a set of non-anticipative strategies for the freezing control to player A, and $\tilde{\Gamma}$ is a set of non-anticipative strategies for player A to player B and the freezing control:
\begin{align}
    &\textrm{N}_A(t)\coloneqq \{ \nu_A :\mathcal{A}(t)\rightarrow \mathcal{M} (t) ~| ~ \forall s\in[t,\tau], \alpha,\bar{\alpha}\in\mathcal{A}(t),\notag\\
    &~~~\text{if } \alpha(\tau)=\bar{\alpha}(\tau) \text{ a.e. } \tau\in[t,s], \notag\\%\text{then } \delta[\alpha](\tau)=\delta[\bar{\alpha}](\tau),\notag\\
    % & \nu_A[\alpha](\tau)=\nu_A[\bar{\alpha}](\tau) \text{ a.e. } \tau\in[t,s] \}.\\
    &~~~\text{then } \nu_A[\alpha](\tau)=\nu_A[\bar{\alpha}](\tau)  \text{ a.e. } \tau\in[t,s] \},
\label{eq:FreezingStrategySet_A}\\
% \end{align}
% Where player B first plays and and the non-anticipative strategy for player A to both player B and the freezing control
% \begin{align}
    &\tilde{\Gamma}(t)\coloneqq \{\tilde\gamma:\mathcal{B}(t)\times\mathcal{M}(t)\rightarrow\mathcal{A}(t) ~| ~ \forall s\in[t,\tau], \beta,\bar{\beta}\in \mathcal{B}(t), \notag\\
    &~~~\mu,\bar{\mu}\in{\mathcal{M}}(t), \text{if } \beta(\tau)=\bar{\beta}(\tau),\mu(\tau)=\bar{\mu}(\tau) \text{ a.e. } \tau\in[t,s], \notag\\
    &~~~\text{then } \tilde\gamma[\beta,\mu](\tau)=\tilde\gamma[\bar{\beta},\bar{\mu}](\tau)  \text{ a.e. } \tau\in[t,s] \}.
\label{eq:playerA_strategySet_aug}
\end{align}

After introducing the auxiliary variable $z\in\R$, define cost $\tilde{J}$ by combining the cost and the constraint of $\tilde\vartheta_1^\pm$: 
\begin{align}
\begin{split}
    &\tilde{J}(t,x,  z,\alpha,\beta,\mu) \coloneqq \max\big\{
    \max_{s\in[t,T]}c(\mathrm{x}(s)),  g(\mathrm{x}(T))-\mathrm{z}(T)
    \big\},
    \end{split}
    \label{eq:ctrl_costDef_TI}
\end{align}
where $(\mathrm{x},\mathrm{z})$ solves, for $s\in[t,T]$,
\begin{align}
\begin{split}
    &\begin{bmatrix}\dot{\mathrm{x}}(s)\\\dot{\mathrm{z}}(s)\end{bmatrix} = \begin{bmatrix}f(\mathrm{x}(s),\alpha(s),\beta(s))\\-L(\mathrm{x}(s),\alpha(s),\beta(s))\end{bmatrix}\mu(s),~
    \begin{bmatrix}\mathrm{x}(t)\\\mathrm{z}(t)\end{bmatrix} =\begin{bmatrix}x\\z\end{bmatrix}.
\end{split}
    \label{eq:dynamics_comb_TI}
\end{align}
Lemma \ref{lemma:equiv_V1} claims that the zero-sum games whose cost is $\tilde J$ are equivalent to $V_1^\pm$ in \eqref{eq:game_V1_upper} and \eqref{eq:game_V1_lower}, which corresponds to $\vartheta_1^\pm$.

\begin{lemma}
    Consider $V_1^\pm$ in \eqref{eq:game_V1_upper} and \eqref{eq:game_V1_lower}, and $\tilde J$ in \eqref{eq:ctrl_costDef_TI}.
    For all $(t,x,z)\in[0,T]\times\R^n\times\R$,
    \begin{align}
        &V_1^+(t,x,z) =    \sup_{\substack{\delta\in\Delta(t),\\\nu_A\in\textrm{N}_A(t)}}\inf_{\alpha\in\mathcal{A}(t)} \tilde{J}(t,x,z,\alpha,\delta[\alpha],\nu_A[\alpha]),
        \label{eq:game_V1_upper_finiteHorizon}\\
    % \end{align}
    % and 
    % \begin{align}
        &V_1^-(t,x,z) =    \inf_{\tilde\gamma\in\tilde{\Gamma}(t)}\sup_{\substack{\beta\in\mathcal{B}(t),\\\mu\in\mathcal{M}(t)}} \tilde{J}(t,x,z,\tilde\gamma[\beta,\mu],\beta,\mu).
        \label{eq:game_V1_lower_finiteHorizon}
    \end{align}
    \label{lemma:equiv_V1}
\end{lemma}
\textbf{Proof.} See Appendix \ref{appen:lemma_equiv_V1}.

\begin{corollary}\textnormal{\textbf{(Equivalent fixed-terminal-time game to the time-invariant Problem \ref{prob:1}})}
    \begin{align}
        \vartheta_1^\pm \equiv \tilde\vartheta_1^\pm,
        \label{eq:equiv_FintieHorizonProblem_game}
    \end{align}
    where $\vartheta_1^+$ is \eqref{eq:def_vartheta1_upper} subject to \eqref{eq:def_vartheta_const1_upper}, $\vartheta_1^-$ is \eqref{eq:def_vartheta1_lower} subject to \eqref{eq:def_vartheta_const1_lower}, $\tilde\vartheta_1^+$ is \eqref{eq:vartheta1_upper_opt_cost} subject to \eqref{eq:vartheta1_upper_opt_const}, and $\tilde\vartheta_1^-$ is \eqref{eq:vartheta1_lower_opt_cost} subject to \eqref{eq:vartheta1_lower_opt_const}. 
    \label{corollary:equiv_FintieHorizonProblem_game}
\end{corollary}
% \textbf{Proof.} See Appendix \ref{appen:lemma_equiv_FintieHorizonProblem1}.
\textbf{Proof.} Let the right hand terms in \eqref{eq:game_V1_upper_finiteHorizon} and \eqref{eq:game_V1_lower_finiteHorizon} be denoted as $W_1^\pm$. By Corollary 5.3 in \cite{altarovici2013general}, $\tilde\vartheta_1^\pm (t,x)=\min z$ subject to $W_1^\pm(t,x,z)\leq 0$. This fact and Lemma \ref{lemma:Equiv_twoCost1} allow us to conclude \eqref{eq:equiv_FintieHorizonProblem_game}. \qed

This corollary remarks that the free-terminal-time games ($\vartheta_1^\pm$) can be converted to fixed-terminal-time games ($\tilde\vartheta_1^\pm$), in which only control signals and strategies have to be specified, since the terminal time is fixed.

In Lemma \ref{lemma:equiv_V1}, $V_1^\pm$ is converted to a fixed-terminal-time game, whose corresponding HJ equation has been investigated in \cite{altarovici2013general}.
This allows us to derive an HJ equation for the time-invariant Problem \ref{prob:1} in Theorem \ref{thm:HJeq_Prob1_TI}.

\begin{theorem}\textnormal{\textbf{(HJ equation for Problem \ref{prob:1} (time-invariant version))}}
    Consider Problem \ref{prob:1} for the time-invariant case.
    For all $(t,x,z)\in[0,T]\times\R^n\times\R$, $V_1^\pm$ in \eqref{eq:game_V1_upper} and \eqref{eq:game_V1_lower} is the unique viscosity solution to the HJ equation:
    \begin{align}
    \begin{split}
        \max&\Big\{c(x)-V_1^\pm (t,x,z), \\
        &V_{1,t}^\pm - \min\big\{0,\bar{H}^\pm (x,z,D_x V_1^\pm,D_z V_1^\pm)\big\} \Big\}=0
        % &V_{1,t}^\pm - \bar{H}_1^{\textnormal{TI}\pm} (x,z,D_x V_1^\pm,D_z V_1^\pm)\Big\} \bigg\}=0
    \end{split}
    \label{eq:HJeq1_TI}
    \end{align}
    in $(0,T)\times\R^n\times\R$,
    where 
    % \begin{align}
    %     \bar{H}^{\textnormal{TI}\pm}_1 (x,z,p,q) \coloneqq \min\Big\{0,\bar{H}^\pm  (x,z,p,q) \Big\}
    % \end{align}
    % for $(x,z,p,q)\in\R^n\times\R\times\R^n\times\R$,
    $\bar{H}^+$ and $\bar{H}^-$ are defined in \eqref{eq:def_Hambar_+} and \eqref{eq:def_Hambar_-}, respectively, without the time dependency,
    and
    \begin{align}
        V_1^\pm (T,x,z) = \max\{c(x),g(x)-z\}
    \end{align}
    on $\{t=T\}\times\R^n\times\R$.
\label{thm:HJeq_Prob1_TI}    
\end{theorem}
\textbf{Proof.} See Appendix \ref{appen:thm_HJeq_Prob1_TI}.
\\Note that the Hamiltonian $\bar{H}^{\pm}$ in \eqref{eq:HJeq1_TI} is time-invariant.
% under the time-invariant conditions, thus it does not have the time dependency.

We observe that the right two terms in the HJ equation \eqref{eq:HJeq1} $\max\{g-z-V_1^\pm,V_{1,t}^\pm - \bar{H}^\pm\}$ become $V_{1,t}^\pm -\min\{0, \bar{H}^\pm\}$ in \eqref{eq:HJeq1_TI}. Note that these two terms are not algebraically equal.

%%%%%%%%%%%%%%%%%%%%%%%%%%%%%%%%%%%%%%%%%%%%%%%%%%%%%%%%%%%%%%%%%%%%%%%%%%%%%%%%%%%%%%%%%%%%%
\subsection{HJ equation for Problem \ref{prob:1} (optimal control setting)}
\label{sec:OptCtrlProb1}

In this subsection, we solve Problem \ref{prob:1} in the optimal control problem setting:
for given initial time and state $(t,x)$,
\begin{align}
 \vartheta_1 (t,x) &\coloneqq \inf_{\alpha\in\mathcal{A}(t)} \max_{\tau\in[t,T]} \int_t^\tau L(s,\mathrm{x}(s),\alpha(s))ds + g(\tau,\mathrm{x}(\tau)),  
\label{eq:def_vartheta1_optCtrl}
\\
& \text{subject to } \quad c(s,\mathrm{x}(s)) \leq 0, \quad s\in [t,\tau],
\label{eq:def_vartheta_const1_optCtrl}
\end{align}
where $\mathrm{x}$ solves 
\begin{align}
    \dot{\mathrm{x}}(s) = f(s,\mathrm{x}(s),\alpha(s)), s\in[t,T], \text{  and  } \mathrm{x}(t)=x.
    \label{eq:dynamics_OptCtrl}
\end{align}
Section \ref{sec:HJeq_gameProb1} and \ref{sec:timeInvProb} present the HJ equations for Problem \ref{prob:1} in the zero-sum game setting. 
By removing player B in the zero-sum game, we can get HJ equations for Problem \ref{prob:1} in the optimal control setting.
Thus, Theorem \ref{thm:HJeq_prob1} and \ref{thm:HJeq_Prob1_TI} imply the following remark.

\begin{remark}\textnormal{\textbf{(HJ equation for Problem \ref{prob:1} (optimal control setting))}}
    Let $V_1$ be the unique viscosity solution to the HJ equation:
    \begin{align}
    \begin{split}
        \max\Big\{c(t,x)-&V_1(t,x,z), g(t,x)-z-V_1(t,x,z), \\ 
        &V_{1,t} - \bar{H} (t,x,z,D_x V_1,D_z V_1)\Big\}=0
        \end{split}
        \label{eq:HJeq1_opt}
    \end{align}
    in $(0,T)\times\R^n\times\R$, where $\bar{H}:[0,T]\times\R^n\times\R\times\R^n\times\R\rightarrow\R$ \begin{align}
        \bar{H}(t,x,z,p,q) \coloneqq \max_{a\in A} - p\cdot f(t,x,a) + q L(t,x,a),
        \label{eq:def_Hambar}
    \end{align}
    and
    \begin{align}
        V_1(T,x,z) = \max\{ c(T,x), g(T,x)-z \}
        % \label{eq:V_terminalValue1}
    \end{align}
    on $\{t=T\}\times\R^n\times\R$. 
    Then, 
    \begin{align}
        \vartheta_1(t,x) = \min z \text{ subject to } V_1(t,x,z)\leq 0,
    \end{align}
    where $\vartheta_1$ is \eqref{eq:def_vartheta1_optCtrl} subject to \eqref{eq:def_vartheta_const1_optCtrl}.
    
    If Problem \ref{prob:1} is time-invariant, $V_1$ is the unique viscosity solution to the HJ equation:
    \begin{align}
    \begin{split}
        \max&\Big\{c(x)-V_1 (t,x,z), \\
        &V_{1,t} - \min\big\{0,\bar{H} (x,z,D_x V_1,D_z V_1)\big\}\Big\}=0
    \end{split}
    \label{eq:HJeq1_optCtrl_TI}
    \end{align}
    in $(0,T)\times\R^n\times\R$,
    where $\bar{H}$ is defined in \eqref{eq:def_Hambar} with ignoring the time dependency,
    and
    \begin{align}
        V_1 (T,x,z) = \max\{c(x),g(x)-z\}
    \end{align}
    on $\{t=T\}\times\R^n\times\R$.
\label{remark:HJeq_prob1_opt}
\end{remark}

%%%%%%%%%%%%%%%%%%%%%%%%%%%%%%%%%%%%%%%%%%%%%%%%%%%%%%%%%%%%%%%%%%%%%%%%%%%%%%%%%%%%%%%%%%%%%%%%%%%%%%%%%%%%%%%%%%%%%%%%%%%%%%%%%%%%%%%%%%%%%%%%%%%%%%%%%%%%%%%%%%%%%%%%%%%%%%%%%%%%%%%%%%%%%%%%%%%%%%%%%%%%%%%%%%%%%%%%%%%%%%%%%%%%%%%%%%%%%%%%
\section{Hamilton-Jacobi Equations for Problem \ref{prob:2}}
\label{sec:HJeq_Prob2Big}

Problem \ref{prob:2} in the zero-sum game setting is defined in Section \ref{sec:problemDefinition}.
Section \ref{sec:HJeq_Prob2} and \ref{sec:HJeq_TIProb2} present HJ equations for the time-varying and time-invariant Problem \ref{prob:2}, respectively.
Section \ref{sec:HJeq_Prob2_OptCtrl} presents HJ equations for Problem \ref{prob:2} and the time-invariant Problem \ref{prob:2} in the optimal control setting.

For Problem \ref{prob:2} in the zero-sum game and optimal control settings, the corresponding HJ equations have been presented in the authors' previous work \cite{donggun2020CDC}.
This section first presents this previous work and then proposes HJ equations for the time-invariant version.

%%%%%%%%%%%%%%%%%%%%%%%%%%%%%%%%%%%%%%%%%%%%%%%%%%%%%%%%%%%%%%%%%%%%%%%%%%%%%%%%
\subsection{HJ equation for the time-varying Problem \ref{prob:2}}
\label{sec:HJeq_Prob2}

This subsection provides an HJ formulation for Problem \ref{prob:2}: solve $\vartheta_2^+$ in \eqref{eq:def_vartheta2_upper} subject to \eqref{eq:def_vartheta_const2_upper} and $\vartheta_2^-$ in \eqref{eq:def_vartheta2_lower} subject to \eqref{eq:def_vartheta_const2_lower}.
For $(t,x,z)\in[0,T]\times\R^n\times\R$, define the augmented value functions corresponding to the upper and lower value functions ($\vartheta_2^\pm$):
\begin{align}
    &V_2^+(t,x,z) \coloneqq \sup_{\delta\in\Delta(t)}\inf_{\alpha\in\mathcal{A}(t)} J_2(t,x,a,\alpha,\delta[\alpha]),
    \label{eq:game_V2_upper}\\
% \end{align}
% and 
% \begin{align}
    &V_2^-(t,x,z) \coloneqq \inf_{\gamma\in\Gamma(t)}\sup_{\beta\in\mathcal{B}(t)} J_2(t,x,a,\gamma[\beta],\beta),
    \label{eq:game_V2_lower}
\end{align}
where cost $J_2:(t,x,z,\alpha,\beta)\rightarrow\R$ is defined as follows:
\begin{align}
\begin{split}
    J_2(t,x&,z,\alpha,\beta) \coloneqq \min_{\tau\in[t,T]}\max\bigg\{ \max_{s\in[t,\tau]}c(s,\mathrm{x}(s)),\\
        &\int_t^\tau L(s,\mathrm{x}(s),\alpha(s),\beta(s))ds + g(\tau,\mathrm{x}(\tau))-z\bigg\},
\end{split}
\label{eq:Prob2_costJ2}
\end{align}
where $\mathrm{x}$ solves \eqref{eq:dynamics} for $(\alpha,\beta)$.
\cite{donggun2020CDC} proved that, for all $(t,x)\in[t,T]\times\R^n$,
\begin{align}
    \vartheta_2^\pm(t,x) =\min z \text{ subject to } V_2^\pm(t,x,z)\leq 0,
\end{align}
and $V_2^\pm$ are the unique viscosity solutions to the HJ equations in Theorem \ref{thm:HJeq_prob2}.
\begin{theorem}\textnormal{\textbf{(HJ equation for Problem \ref{prob:2}) \cite{donggun2020CDC}}}
    For all $(t,x,z)\in[0,T]\times\R^n\times\R$, $V_2^\pm$ in \eqref{eq:game_V2_upper} and \eqref{eq:game_V2_lower} are the unique viscosity solutions to the HJ equations:
    \begin{align}
    % \begin{split}
        \max\Big\{c&(t,x)-V_2^\pm(t,x,z), \min\big\{g(t,x)-z-V_2^\pm(t,x,z), \notag\\ 
        &V_{2,t}^\pm - \bar{H}^\pm (t,x,z,D_x V_2^\pm,D_z V_2^\pm)\big\}\Big\}=0
        % \end{split}
        \label{eq:HJeq2}
    \end{align}
    in $(0,T)\times\R^n\times\R$, where $\bar{H}^\pm$ are defined in \eqref{eq:def_Hambar_+} and \eqref{eq:def_Hambar_-},
    and
    \begin{align}
        V_2^\pm(T,x,z) = \max\{ c(T,x), g(T,x)-z \}
        \label{eq:V_terminalValue2}
    \end{align}
    on $\{t=T\}\times\R^n\times\R$.
    \label{thm:HJeq_prob2}
\end{theorem}
We observe that the difference between the two types of HJ equations for $V_1^\pm$ and $V_2^\pm$ is that the minimum operation in \eqref{eq:HJeq2} for $V_2^\pm$ is replaced by the maximum operation in \eqref{eq:HJeq1}. 
This is from the difference between $\vartheta_1^\pm$ and $\vartheta_2^\pm$: $\vartheta_1^\pm$ in \eqref{eq:def_vartheta1_upper} and \eqref{eq:def_vartheta1_lower} have $\max_{\tau}$ operation, and $\vartheta_2^\pm$ in \eqref{eq:def_vartheta2_upper} and \eqref{eq:def_vartheta2_lower} have $\min_{\tau}$ operation.

%%%%%%%%%%%%%%%%%%%%%%%%%%%%%%%%%%%%%%%%%%%%%%%%%%%%%%%%%%%%%%%%%%%%%%%%%%%%%%%%
\subsection{HJ equation for Problem \ref{prob:2} (time-invariant case)}
\label{sec:HJeq_TIProb2}

% Consider time-invariant problems.
Through similar analysis to that in Section \ref{sec:timeInvProb}, this subsection derives the HJ equations for the time-invariant Problem \ref{prob:2}.
For Problem \ref{prob:1}, the freezing control signal $\mu:[t,T]\rightarrow[0,1]$ allows to convert to the fixed-terminal-time problems by replacing the maximum over $\tau$ operation in Problem \ref{prob:1} to the supremum over the freezing control signal or strategy.
Instead, Problem \ref{prob:2} is specified in terms of the minimum over $\tau$ operation, which will be replaced by the infimum over the freezing control signal or strategy.

Consider two fixed-terminal-time problems:
\begin{align}
 \begin{split}
%\sup_{\substack{\delta\in \Delta(t) \\ \nu_A\in \Delta_d(t)}}
    \tilde\vartheta_2^+ (t&,x) \coloneqq \sup_{\tilde\delta\in \tilde\Delta(t)} \inf_{\alpha\in\mathcal{A}(t),\mu\in\mathcal{M}(t)} \\
    &\int_t^T L(\mathrm{x}(s),\alpha(s),\tilde\delta[\alpha,\mu](s)) \mu(s)ds + g(\mathrm{x}(T)),
    \label{eq:vartheta2_upper_opt_cost}
\end{split}\\
&\quad\quad\quad \text{subject to } c(\mathrm{x}(s))\leq 0, s\in[t,T],
\label{eq:vartheta2_upper_opt_const}         
\end{align}
where $\mathrm{x}$ solves \eqref{eq:dynamics_timeInvariant} for $(\alpha,\tilde\delta[\alpha,\mu],\mu)$, $\mathcal{M}$ is defined in \eqref{eq:ctrlSet_freezing}, and $\tilde\delta$ is the non-anticipative strategy for player B to both player A and the freezing control:
\begin{align}
    \tilde{\Delta}&(t)\coloneqq \{\tilde\delta:\mathcal{A}(t)\times\mathcal{M}(t)\rightarrow\mathcal{B}(t) ~| ~ \forall s\in[t,\tau], \alpha,\bar{\alpha}\in \mathcal{A}(t), \notag\\
    &\mu,\bar{\mu}\in \mathcal{M}(t), \text{if } \alpha(\tau)=\bar{\alpha}(\tau),\mu(\tau)=\bar{\mu}(\tau) \text{ a.e. } \tau\in[t,s], \notag\\
    &\text{then } \tilde\delta[\beta,\mu](\tau)=\tilde\delta[\bar{\beta},\bar{\mu}](\tau)  \text{ a.e. } \tau\in[t,s] \};
\label{eq:playerB_strategySet_aug}
\end{align}
\begin{align}
\begin{split}
    \tilde\vartheta_1^- (t&,x) \coloneqq\inf_{\gamma\in \Gamma(t),\nu_B\in \textrm{N}_B(t)} \sup_{\beta\in\mathcal{B}(t)} \\
    &\int_t^T L(\mathrm{x}(s),\gamma[\beta](s),\beta(s)) \nu_B[\beta](s)ds + g(\mathrm{x}(T)),
    \label{eq:vartheta2_lower_opt_cost}
\end{split}\\
&\quad\quad\quad \text{subject to } c(\mathrm{x}(s))\leq 0, s\in[t,T],
\label{eq:vartheta2_lower_opt_const}
\end{align}
where $\mathrm{x}$ solves \eqref{eq:dynamics_timeInvariant} for $(\gamma[\beta],\beta,\nu_B[\beta])$, and the non-anticipative strategy for the freezing control to player B is 
\begin{align}
    \textrm{N}_B(t)&\coloneqq \{ \nu_B :\mathcal{B}(t)\rightarrow \mathcal{M} (t) ~| ~ \forall s\in[t,\tau], \beta,\bar{\beta}\in\mathcal{B}(t),\notag\\
    &\text{if } \beta(\tau)=\bar{\beta}(\tau) \text{ a.e. } \tau\in[t,s], \notag\\%\text{then } \delta[\alpha](\tau)=\delta[\bar{\alpha}](\tau),\notag\\
    % & \nu_A[\alpha](\tau)=\nu_A[\bar{\alpha}](\tau) \text{ a.e. } \tau\in[t,s] \}.\\
    &\text{then } \nu_B[\beta](\tau)=\nu_B[\bar{\beta}](\tau)  \text{ a.e. } \tau\in[t,s] \}.
\label{eq:FreezingStrategySet_B}
\end{align}

Recall $\tilde J$ in \eqref{eq:ctrl_costDef_TI}, which contains the cost and the constraint of $\tilde{\vartheta}_2^\pm$. 
Following similar steps of the proof of Lemma \ref{lemma:equiv_V1}, Lemma \ref{lemma:equiv_V2} can be proved.
\begin{lemma}
    Recall $\tilde J$ in \eqref{eq:ctrl_costDef_TI}, and consider $V_2^\pm$ in \eqref{eq:game_V2_upper} and \eqref{eq:game_V2_lower}.
    For all $(t,x,z)\in[0,T]\times\R^n\times\R$,
    \begin{align}
        &V_2^+(t,x,z) =    \sup_{\tilde\delta\in\tilde\Delta(t)}\inf_{\substack{\alpha\in\mathcal{A}(t),\\\mu\in\mathcal{M}(t)}} \tilde{J}(t,x,z,\alpha,\tilde\delta[\alpha,\mu],\mu),
        \label{eq:game_V2_upper_finiteHorizon}\\
    % \end{align}
    % and 
    % \begin{align}
        &V_2^-(t,x,z) = \inf_{\substack{\gamma\in\Gamma(t),\\\nu_B\in\textrm{N}_B(t)}} \sup_{\beta\in\mathcal{B}(t)} \tilde{J}(t,x,z,\gamma[\beta],\beta,\nu_B[\beta]).
        \label{eq:game_V2_lower_finiteHorizon}
    \end{align}
    \label{lemma:equiv_V2}
\end{lemma}
% \textbf{Proof.} See Appendix \ref{appen:lemma_equiv_V2}.

By combining the HJ formulation for the fixed-terminal-time problems \cite{altarovici2013general} and Lemma \ref{lemma:equiv_V2}, the HJ equation for the time-invariant Problem \ref{prob:2} is derived in Theorem \ref{thm:HJeq_Prob2_TI}.
The proof for Theorem \ref{thm:HJeq_Prob2_TI} is analogous to the proof for Theorem \ref{thm:HJeq_Prob1_TI}. 

\begin{theorem}\textnormal{\textbf{(HJ equation for Problem \ref{prob:2} (time-invariant version))}}
    Consider Problem \ref{prob:2} in the time-invariant case.
    For all $(t,x,z)\in[0,T]\times\R^n\times\R$, $V_2^\pm$ in \eqref{eq:game_V2_upper} and \eqref{eq:game_V2_lower} are the unique viscosity solutions to the HJ equations:
    \begin{align}
    \begin{split}
        \max&\Big\{c(x)-V_2^\pm (t,x,z), \\
        &V_{2,t}^\pm - \max\big\{0,  \bar{H}^\pm (x,z,D_x V_2^\pm,D_z V_2^\pm)\big\}\Big\}=0
    \end{split}
    \label{eq:HJeq2_TI}
    \end{align}
    in $(0,T)\times\R^n\times\R$,
    where $\bar{H}^+$ and $\bar{H}^-$ are defined in \eqref{eq:def_Hambar_+} and \eqref{eq:def_Hambar_-}, respectively, without the time dependency,
    and
    \begin{align}
        V_2^\pm (T,x,z) = \max\{c(x),g(x)-z\}
    \end{align}
    on $\{t=T\}\times\R^n\times\R$.
\label{thm:HJeq_Prob2_TI}    
\end{theorem}
% \textbf{Proof.} See Appendix \ref{appen:thm_HJeq_Prob2_TI}.

In comparison between the HJ equations for Problem \ref{prob:2} and its time-invariant version, $\min\{g-z-V_2^\pm,V_{2,t}^\pm -\bar{H}^\pm\}$ in \eqref{eq:HJeq2} becomes $V_{2,t}^\pm - \max\{0, \bar{H}^\pm\}$ in \eqref{eq:HJeq2_TI}.
Note that the difference between Problem \ref{prob:1} and \ref{prob:2} leads to the difference in HJ equations for the time-invariant problems: \eqref{eq:HJeq1_TI} has the term $V_{1,t}^\pm - \min\{0, \bar{H}^\pm\}$, but \eqref{eq:HJeq2_TI} has the term $V_{2,t}^\pm - \max\{0, \bar{H}^\pm\}$.

%%%%%%%%%%%%%%%%%%%%%%%%%%%%%%%%%%%%%%%%%%%%%%%%%%%%%%%%%%%%%%%%%%%%%%%%%%%%%%%%
\subsection{HJ equation for Problem \ref{prob:2} (optimal control setting)}
\label{sec:HJeq_Prob2_OptCtrl}

In this subsection, we solve Problem \ref{prob:2} in the optimal control setting:
for given initial time and state $(t,x)$,
\begin{align}
 \vartheta_2 (t,x) &\coloneqq \inf_{\alpha\in\mathcal{A}(t)} \min_{\tau\in[t,T]} \int_t^\tau L(s,\mathrm{x}(s),\alpha(s))ds + g(\tau,\mathrm{x}(\tau)),  
\label{eq:def_vartheta2_optCtrl}
\\
& \text{subject to } \quad c(s,\mathrm{x}(s)) \leq 0, \quad s\in [t,\tau],
\label{eq:def_vartheta_const2_optCtrl}
\end{align}
where $\mathrm{x}$ solves \eqref{eq:dynamics_OptCtrl}.
By removing the contribution of player B in Theorem \ref{thm:HJeq_prob2} and \ref{thm:HJeq_Prob2_TI}, we solve $\vartheta_2$ using the HJ equations in the following remark.

\begin{remark}\textnormal{\textbf{(HJ equation for Problem \ref{prob:2} (optimal control setting))}}
    Let $V_2$ be the unique viscosity solution to the HJ equation \cite{donggun2020CDC}:
    \begin{align}
    % \begin{split}
        \max\Big\{&c(t,x)-V_2(t,x,z), \min\big\{ g(t,x)-z-V_2(t,x,z), \notag\\ 
        &V_{2,t} - \bar{H} (t,x,z,D_x V_2,D_z V_2)\big\} \Big\}=0
        % \end{split}
        \label{eq:HJeq2_opt}
    \end{align}
    in $(0,T)\times\R^n\times\R$, where $\bar H$ is defined in \eqref{eq:def_Hambar},
    and
    \begin{align}
        V_2(T,x,z) = \max\{ c(T,x), g(T,x)-z \}
        % \label{eq:V_terminalValue2}
    \end{align}
    on $\{t=T\}\times\R^n\times\R$. 
    Then,
    \begin{align}
        \vartheta_2(t,x) = \min z \text{ subject to }V_2(t,x,z)\leq 0,
    \end{align}
    where $\vartheta_2$ is \eqref{eq:def_vartheta2_optCtrl} subject to \eqref{eq:def_vartheta_const2_optCtrl}.
    
    If Problem \ref{prob:2} is time-invariant, $V_2$ is the unique viscosity solution to the HJ equation:
    \begin{align}
    \begin{split}
        \max&\Big\{c(x)-V_2 (t,x,z), \\
        &V_{2,t} - \max\big\{0,\bar{H} (x,z,D_x V_2,D_z V_2)\big\}\Big\}=0
    \end{split}
    \label{eq:HJeq2_optCtrl_TI}
    \end{align}
    in $(0,T)\times\R^n\times\R$,
    where $\bar{H}$ is defined in \eqref{eq:def_Hambar} without the time dependency,
    and
    \begin{align}
        V_2 (T,x,z) = \max\{c(x),g(x)-z\}
    \end{align}
    on $\{t=T\}\times\R^n\times\R$.
\label{remark:HJeq_prob2_opt}
\end{remark}

%%%%%%%%%%%%%%%%%%%%%%%%%%%%%%%%%%%%%%%%%%%%%%%%%%%%%%%%%%%%%%%%%%%%%%%%%%%%%%%%%%%%%%%%%%%%%%%%%%%%%%%%%%%%%%%%%%%%%%%%%%%%%%%%%%%%%%%%%%%%%%%%%%%%%%%%%%%%%%%%
%%%%%%%%%%%%%%%%%%%%%%%%%%%%%%%%%%%%%%%%%%%%%%%%%%%%%%%%%%%%%%%%%%%%%%%%%%%%%%%%
\section{Optimal Control Signal and Strategy}
\label{sec:optCtrl_strategy}

The optimal control signal or strategy for Problems \ref{prob:1} and \ref{prob:2} are specified by the HJ equations in Section \ref{sec:HJeq_game} and \ref{sec:HJeq_Prob2Big}.
This section utilizes the HJ equations in Theorem \ref{thm:HJeq_prob1} and \ref{thm:HJeq_prob2}, and the method in this section can be simply extended for the other HJ equations in Theorems \ref{thm:HJeq_Prob1_TI} and \ref{thm:HJeq_Prob2_TI}, and Remarks \ref{remark:HJeq_prob1_opt} and \ref{remark:HJeq_prob2_opt}.

Recall $V_i^\pm$ ($i=1,2$) defined in \eqref{eq:game_V1_upper}, \eqref{eq:game_V1_lower}, \eqref{eq:game_V2_upper}, \eqref{eq:game_V2_lower}, and suppose $V_i^\pm$ is computed from the HJ equations in Theorems \ref{thm:HJeq_prob1} and \ref{thm:HJeq_prob2}.

Lemmas \ref{lemma:equiv_V1} and \ref{lemma:equiv_V2} imply the following remark.
\begin{remark}[Find $\vartheta_i^\pm$ from $V_i^\pm$] For initial time $t=0$ and state $x\in\R^n$, 
    \begin{align}
        (\mathrm{x}_*(0),\mathrm{z}_*(0)) = (x, \vartheta_i^\pm(0,x) ) ,
    \end{align}
    where $(\mathrm{x}_*,\mathrm{z}_*)$ is an optimal trajectory for $V_i^\pm$.
\label{remark:Optctrl_vartheta}
\end{remark}

% Remark \ref{remark:Optctrl_vartheta} presents optimal control signal and strategy for the games $\vartheta_i^\pm$ using those for $V_i^\pm$.
% In the rest of this section, we analyze the optimal control and strategy for $V_i^\pm$.

With the initial augmented state $(\mathrm{x}_*(0),\mathrm{z}_*(0))$, the optimal control and strategy can be found at $(t,\mathrm{x}_*(t),\mathrm{z}_*(t))$, and the optimal state trajectory is also updated by solving the ODE \eqref{eq:dynamics_comb}.

Define $\tilde{H}_i^\pm: A\times B\rightarrow \R$ for a fixed $(t,x,z)\in (0,T)\times\R^n\times\R$
\begin{align}
\begin{split}
    \tilde{H}_i^\pm(a,b) \coloneqq &-D_x V_i^\pm(t,x,z) \cdot f(t,x,a,b) \\
    &+ D_z V_i^\pm(t,x,z) L(t,x,a,b)
\end{split}
\label{eq:tilde_Hamiltonian}
\end{align}
for $i=1,2$, thus 
\begin{align}
    \bar{H}^+(t,x,z,D_x V_i^+,D_z V_i^+) = \max_{a\in A} \min_{b\in B} \tilde{H}_i^+(a,b),\\
    \bar{H}^- (t,x,z,D_x V_i^-,D_z V_i^-) = \min_{b\in B} \max_{a\in A} \tilde{H}_i^- (a,b),
\end{align}
where $\bar{H}^+$ and $\bar{H}^-$ are defined in \eqref{eq:def_Hambar_+} and \eqref{eq:def_Hambar_-}, respectively.
% At each $(t,x,z)$, the minimax solution $(a_*,b_*)$ for the Hamiltonian $\bar{H}_V^\pm$ is not only optimal, but any pair of $(a,b)$ satisfying the HJ equations in \eqref{eq:HJeq1} or \eqref{eq:HJeq2} is also optimal.
In this section, we omit $(t,x,z)$ to simplify notation.
Using the notation with $\tilde{H}_i^\pm$ \eqref{eq:tilde_Hamiltonian}, the HJ equation \eqref{eq:HJeq1} for $V_1^+$ is equal to
\begin{align}
    \max\{c-V_1^+, g-z-V_1^+,V_{1,t}^+ -\max_{a\in A}\min_{b\in B}\tilde{H}_1^+ (a,b) \}=0.
    \label{eq:HJeq1+1}
\end{align}

The HJ equation implies that optimal control signal or strategy is determined by the gradient information at the current time and state $(t,x,z)$, but the past history of the state trajectory and optimal control signals is not necessary.
For example, in $V_1^+$, $\alpha_*(t) = a_*$, $\delta_*[\alpha_*](t)= b_*$ where $a_*$ and $b_*$ are solutions to
\begin{align}
    & \max\{c-V_1^+, g-z-V_1^+,V_{1,t}^+ -\tilde{H}_1^+ (a,b) \}=0.
    \label{eq:HJeq1+2}
\end{align}
at $(t,x,z)$. In other words, it is sufficient to specify optimal controls for player A and B in $(0,T)\times\R^n\times\R$ to generate the optimal control signal or strategy.
% Since the optimal controls for player A and B in $(0,T)\times\R^n\times\R$ generate the optimal control signal or strategy, it is sufficient to find optimal controls for player A and B in this section.
The maxmini solution $(a_*,b_*)$ for the Hamiltonian $\bar{H}_V^+$ ($\max_{a\in A}\min_{b\in B}\tilde{H}_1^+ (a,b)$) is certainly optimal, but there are more solutions.
Similarly, for $V_1^-$ or $V_2^\pm$, any pair of $(a_*,b_*)$ satisfying the corresponding HJ equation is optimal.

% At each $(t,x,z)\in(0,T)\times \R^n\times\R$, a control or strategy is optimal if it achieves the cost-to-go value function $V_i^\pm$.
% To explain this,
% , but any pair of $(a,b)$ satisfying
% \begin{align}
%     \max\{c-V_1^\pm, g-z-V_1^\pm,V_{1,t}^\pm -\tilde{H}_1^\pm (a,b) \}=0
%     \label{eq:HJeq1+2}
% \end{align}
% is also optimal, where player A plays first and player B plays second.

% satisfying \eqref{eq:HJeq1+2} is also optimal, where player B plays first and player A plays second.

% For $V_2^\pm$ in \eqref{eq:game_V2_upper} and \eqref{eq:game_V2_lower}, any pair of $(a,b)$ satisfying 
% \begin{align}
%     \max\{c-V_2^\pm, \min\{ g-z-V_2^\pm,V_{2,t}^\pm -\tilde{H}_2^\pm (a,b)\} \}=0
%     \label{eq:HJeq2+1}
% \end{align}
% is optimal, where player A plays first for $V_2^+$ and player B plays first for $V_2^-$.

In the HJ equation \eqref{eq:HJeq1+1} (or \eqref{eq:HJeq1}) for $V_1^\pm$, we have three terms: $c-V_1^\pm$, $g-z-V_1^\pm$, or $V_{1,t}^\pm-\bar{H}^\pm$, and, at least, one of these terms is zero.
% In the HJ equation \eqref{eq:HJeq2+1} (or \eqref{eq:HJeq2}) for $V_2^\pm$, we also have the same three terms where $V_1^\pm$ is replaced by $V_2^\pm$.
By considering which term is bigger or smaller among the three terms, all possible optimal controls for $\vartheta_1^\pm$ ($V_1^\pm$) and $\vartheta_2^\pm$ ($V_2^\pm$) is derived in Remark \ref{remark:OptCtrl_Strat}.

Although $(a_*,b_*)$ satisfying the HJ equation is optimal, we need to consider the order of players: player A plays first in $V_i^+$ but player B plays first in $V_i^-$.
In $V_i^+$, we first find a set of optimal control for player A, and then investigate a set of optimal control for player B when player A applies its optimal control.
% if $0=c-V_1^+ \geq g-z-V_1^+ \geq V_{1,t}^+ - \bar{H}^+$ at $(t,x,z)$, any player A's control satisfying \eqref{eq:V1+_case1_a} is optimal, and then we investigate a set of optimal control for player B when player A applies its optimal control: all $b_*\in B$ is optimal.
% In this case, all $(a,b)$ pairs satisfying both \eqref{eq:V1+_case1_a} and \eqref{eq:V1+_case1_b} satisfy \eqref{eq:HJeq1+2}.
On the other hand, in $V_1^-$, 
% if $0=c-V_1^- \geq g-z-V_1^- \geq V_{1,t}^- - \bar{H}^-$, 
we first investigate a set of optimal control for player B, and then find a set of optimal control for player A after applying an optimal control of player B.
% : the optimality conditions are \eqref{eq:V1-_case1_b} and \eqref{eq:V1-_case1_a}, which satisfies \eqref{eq:HJeq1+2}.
Based on this argument, Remark \ref{remark:OptCtrl_Strat} presents optimal controls for $V_i^+$ according to classification, and optimal controls for $V_i^-$ can be analogously extended.

\begin{remark}[Optimal controls for $V_i^+$]
    Fix $(t,x,z)\in(0,T)\times\R^n\times\R$.
    
Optimal controls for $\vartheta_1^+$ ($V_1^+$) are the following:
\begin{enumerate}
    \item Case 1: $\max\{c-V_1^+ , g-z- V_1^+\}\geq V_{1,t}^+ -\bar{H} ^+ $
    \begin{align}
        &a_* \in \{a\in A~|~ V_{1,t}^+ -\min_{b\in B} \tilde{H}_1^+(a,b) \leq 0 \},\label{eq:V1+_case1_a}\\
        &b_* \in B;\label{eq:V1+_case1_b}
    \end{align}
    
    \item Case 2: $\max\{c-V_1^+ , g-z- V_1^+\} < V_{1,t}^+ -\bar{H} ^+$
    \begin{align}
        &a_* \in \arg\max_{a\in A}\min_{b\in B}\tilde{H}_1^+ (a,b),\label{eq:V1+_case2_a}\\
        &b_* \in \arg\min_{b\in B} \tilde{H}_1^+ (a_*,b).\label{eq:V1+_case2_b}
    \end{align}
\end{enumerate}
    
% Optimal controls for $\vartheta_1^-$ ($V_1^-$) are the following:
% \begin{enumerate}
%     \item Case 1: $\max\{c-V_1^- , g-z- V_1^-\}\geq V_{1,t}^- -\bar{H} ^-$
%     \begin{align}
%         &b_* \in B,\label{eq:V1-_case1_b}\\
%         &a_* \in \{a\in A~|~ V_{1,t}^- - \tilde{H}_1^- (a,b_*) \leq 0 \};\label{eq:V1-_case1_a}
%     \end{align}
    
%     \item Case 2: $\max\{c-V_1^- , g-z- V_1^-\}<V_{1,t}^- -\bar{H} ^-$
%     \begin{align}
%         &b_* \in \arg\min_{b\in B}\max_{a\in A}\tilde{H}_1^- (a,b),\label{eq:V1-_case2_b}\\
%         &a_* \in \arg\max_{a\in A} \tilde{H}_1^- (a , b_*).\label{eq:V1-_case2_a}
%     \end{align}
% \end{enumerate}
    
Optimal controls for $\vartheta_2^+$ $(V_2^+)$ are the following:
\begin{enumerate}
    \item Case 1: $c-V_2^+ \geq V_{2,t}^+ -\bar{H} ^+$
    
    \eqref{eq:V1+_case1_a} and \eqref{eq:V1+_case1_b} where $V_1^+$ and $\tilde{H}_1^+$ are replaced by $V_2^+$ and $\tilde{H}_2^+$, respectively;
    
    \item Case 2: $g-z-V_2^+ \geq V_{2,t}^+ -\bar{H} ^+ \geq c-V_2^+$
    
    \eqref{eq:V1+_case2_a} and \eqref{eq:V1+_case2_b} where $V_1^+$ and $\tilde{H}_1^+$ are replaced by $V_2^+$ and $\tilde{H}_2^+$, respectively;
    
    \item Case 3: $V_{2,t}^+ -\bar{H} ^+ \geq \max\{c-V_2^+, g-z-V_2^+\}$
    \begin{align}
        &a_* \in A, \label{eq:V2+_case3_a}\\
        &b_* \in \{b\in B~|~ V_{2,t}^+ - \tilde{H}_2^+ (a_*,b) \geq 0 \}. \label{eq:V2+_case3_b}
    \end{align}
\end{enumerate}

% Optimal controls for $\vartheta_2^-$ $(V_2^-)$ is the following:
% \begin{enumerate}
%     \item Case 1: $c-V_2^- \geq V_{2,t}^- -\bar{H} ^-$
    
%     \eqref{eq:V1-_case1_b} and \eqref{eq:V1-_case1_a} where $V_1^-$ and $\tilde{H}_1^-$ are replaced by $V_2^-$ and $\tilde{H}_2^-$, respectively;
    
%     \item Case 2: $g-z-V_2^- \geq V_{2,t}^- -\bar{H} ^- \geq c-V_2^-$
    
%     \eqref{eq:V1-_case2_a} and \eqref{eq:V1-_case2_b} where $V_1^-$ and $\tilde{H}_1^-$ are replaced by $V_2^-$ and $\tilde{H}_2^-$, respectively;
    
%     \item Case 3: $V_{2,t}^- -\bar{H} ^- \geq \max\{c-V_2^-, g-z-V_2^-\}$
%     \begin{align}
%         &b_* \in \{b\in B~|~ V_{2,t}^- - \max_{a\in A}\tilde{H}_2^- (a,b) \geq 0 \}, \label{eq:V2-_case3_b}\\
%         &a_* \in A. \label{eq:V2-_case3_a}
%     \end{align}
% \end{enumerate}

\label{remark:OptCtrl_Strat}
\end{remark}

%%%%%%%%%%%%%%%%%%%%%%%%%%%%%%%%%%%%%%%%%%%%%%%%%%%%%%%%%%%%%%%%%%%%%%%%%%%%%%%%
\section{Numerical Computation for the Hamilton-Jacobi equation}
\label{sec:NumericalAlgorithm}

In this section, we present a numerical algorithm based on level set methods \cite{mitchell2005toolbox} to compute the solutions to the four HJ equations for Problems \ref{prob:1} and \ref{prob:2}.
Algorithm \ref{alg:HJ_computation_prob1} deals with the HJ equations for Problem \ref{prob:1} (\eqref{eq:HJeq1}, \eqref{eq:HJeq1_TI}, \eqref{eq:HJeq1_opt}, \eqref{eq:HJeq1_optCtrl_TI}), and Algorithm \ref{alg:HJ_computation_prob2} deals with the HJ equations for Problem \ref{prob:2} (\eqref{eq:HJeq2}, \eqref{eq:HJeq2_TI}, \eqref{eq:HJeq2_opt}, \eqref{eq:HJeq2_optCtrl_TI}).
Level set methods have been utilized to solve a variety of HJ formulations \cite{Mitchell03,fisac2015reach,donggun2020CDC}.

Algorithm \ref{alg:HJ_computation_prob1} solves the HJ equation \eqref{eq:HJeq1} in two steps. At line 6 in Algorithm \ref{alg:HJ_computation_prob1}, we first compute the HJ PDE ($V_{1,t}^\pm - \bar{H} ^\pm(t,x,z,D_x V_1^\pm,D_z V_1^\pm)=0$), and line 7 in Algorithm \ref{alg:HJ_computation_prob1} replaces $V_1^\pm$ by one of $c(t,x)$, $g(t,x)-z$ or itself to satisfy the HJ equation \eqref{eq:HJeq1}.

For solving the HJ PDE at step 1, the Lax-Friedrichs scheme \cite{crandall1984two} is utilized on the temporal discretization $\{t_0=0,...,t_K=T\}$ and the spatial discretization $\{(x_0,z_0),...,(x_N,z_N)\}\subset \R^n\times\R$:
\begin{align}
    V_1^\pm(t_k,x_i,z_i) = V_1^\pm (t_{k+1},x_i,z_i) -\Delta_k \hat{\bar{H}}^\pm (\phi_x^+,\phi_x^-,\phi_z^+,\phi_z^-),
    \label{eq:first_order_LaxFriedrichs}
\end{align}
where $\Delta_k=t_{k+1}-t_k$, $(\phi_x^\pm,\phi_z^\pm)$ are numerical approximation for $(D_x V_1^\pm, D_z V_1^\pm)$ (gradients with respect to $(x,z)$ at $(t_{k+1},x_i,z_i)$), and
\begin{align}
    \hat{\bar{H}}^\pm (\phi_x^+,\phi_x^-,\phi_z^+&,\phi_z^-) = \bar{H}^\pm\big(t_{k+1},x_i,z_i,\frac{\phi_x^+ + \phi_x^-}{2},\frac{\phi_z^+ + \phi_z^-}{2}\big) \notag\\
    &-\alpha_x \cdot \frac{\phi_x^+ + \phi_x^-}{2}-\alpha_z  \frac{\phi_z^+ + \phi_z^-}{2},
\end{align}
where $\hat {H}^\pm$ are defined in \eqref{eq:def_Hambar_+} and \eqref{eq:def_Hambar_-}, and $\alpha_x=(\alpha_{x_1},...,\alpha_{x_n})$ ($\alpha_{x_i} =\max |D_{p_i} \bar{H} ^\pm |$) and $\alpha_z (=\max | D_q \bar{H} ^\pm |)$ are dissipation coefficients for numerical viscosity, based on the partial derivatives of $\bar{H}^\pm$ \cite{Osher02}.
The fifth-order accurate HJ WENO (weighted essentially nonoscillatory) method \cite{Osher02} is used for the gradient $\phi_x^\pm , \phi_z^\pm$.
In \eqref{eq:first_order_LaxFriedrichs}, the first-order Euler method is used for the temporal partial derivative, but higher-order methods, such as third-order accurate TVD (total variation diminishing) Runge-Kutta (RK) \cite{shu1988efficient} can be used.
\cite{Osher02} provided the empirical observation that level set methods are sensitive to spatial accuracy, thus high-order scheme for spatial derivatives is desired, but high-order approximation for temporal derivatives does not significantly increase the accuracy.

% At step 2, line 7 in Algorithm \ref{alg:HJ_computation_prob1} updates $V_1^\pm$ to the maximum of $c, g-z$ and itself to satisfy the HJ equation \eqref{eq:HJeq1}.

For the time-invariant Problem \ref{prob:1}, line 9 in Algorithm \ref{alg:HJ_computation_prob1} solves the HJ equation \eqref{eq:HJeq1_TI} whose Hamiltonian has the minimum with 0 operation.
Then, line 10 in Algorithm \ref{alg:HJ_computation_prob1} updates $V_1^\pm$ with the maximum between $c$ and itself to satisfy the HJ equation \eqref{eq:HJeq1_TI} without considering $g-z$ term.

For Problem \ref{prob:1} in optimal control setting, Algorithm \ref{alg:HJ_computation_prob1} also works with utilizing the correct Hamiltonian $\bar{H}$ \eqref{eq:def_Hambar} instead of $\bar{H}^\pm$ (\eqref{eq:def_Hambar_+} and \eqref{eq:def_Hambar_-}).

\begin{algorithm}
\caption{Computing the solution $V_1^\pm$ or $V_1$ to the HJ equations for Problem \ref{prob:1} in the zero-sum game and optimal control settings. This algorithm deals with the four HJ equations: \eqref{eq:HJeq1}, \eqref{eq:HJeq1_TI}, \eqref{eq:HJeq1_opt}, and \eqref{eq:HJeq1_optCtrl_TI}.}
\begin{algorithmic}[1]
\State \textbf{Input:} {the temporal discretization: $\{t_0=0,...,t_K=T\}$, the spatial discretization: $\{(x_0,z_0),...,(x_N,z_N)\}$}
\State \textbf{Output:} {$V_1^\pm$ (or $V_1$)}
\State $V_1^\pm (\text{or }V_1) (T,x_i,z_i) \leftarrow \max\{c(T,x_i),g(T,x_i)-z_i\},\forall i$
\For {$k\in\{K-1,...,0\}$}
    \Case{solving the HJ equations \eqref{eq:HJeq1} or \eqref{eq:HJeq1_opt}}
        \State \Longunderstack[l]{$V_1^\pm (\text{or }V_1)(t_{k},x_i,z_i) \leftarrow V_1^\pm(\text{or }V_1) (t_{k+1},x_i,z_i)-$\\
        \quad\quad\quad $\Delta_k \hat{\bar{H}}^\pm (\text{or }\hat{\bar{H}}) (\phi_x^+,\phi_x^-,\phi_z^+,\phi_z^-), \forall i$}
        
        \State \Longunderstack[l]{$V_1^\pm (\text{or }V_1)(t_{k},x_i,z_i) \leftarrow \max\{c(t_k, x_i),$\\
        \quad\quad\quad $g(t_k,x_i)-z_i,V_1^\pm(\text{or }V_1) (t_{k},x_i,z_i)\}, \forall i$}
    \EndCase
    \Case{solving the HJ equations \eqref{eq:HJeq1_TI} or \eqref{eq:HJeq1_optCtrl_TI}}
        \State \Longunderstack[l]{$V_1^\pm(\text{or }V_1) (t_{k},x_i,z_i) \leftarrow V_1^\pm(\text{or }V_1) (t_{k+1},x_i,z_i)-$\\
        \quad\quad $\Delta_k \min\{0,\hat{\bar{H}}^\pm(\text{or }\hat{\bar{H}})  (\phi_x^+,\phi_x^-,\phi_z^+,\phi_z^-)\}, \forall i$}
        
        \State \Longunderstack[l]{$V_1^\pm (\text{or }V_1)(t_{k},x_i,z_i) \leftarrow \max\{c(t_k, x_i),$\\
        \quad\quad\quad $V_1^\pm(\text{or }V_1) (t_{k},x_i,z_i)\}, \forall i$}
    \EndCase
\EndFor
\end{algorithmic}
\label{alg:HJ_computation_prob1}
\end{algorithm}

Algorithm \ref{alg:HJ_computation_prob2} numerically solves Problem \ref{prob:2} in zero-sum game and optimal control setting.
At line 6 in Algorithm \ref{alg:HJ_computation_prob2}, we compute the HJ PDE ($V_{2,t}^\pm - \bar{H} ^\pm(t,x,z,D_x V_2^\pm,D_z V_2^\pm)=0$) in the HJ equation \eqref{eq:HJeq2}.
Line 7 and 8 in Algorithm \ref{alg:HJ_computation_prob2} update $V_2^\pm$ by choosing among $c,g-z$ and itself so that the HJ equation \eqref{eq:HJeq2} holds.

For the time-invariant Problem \ref{prob:2}, line 10 in Algorithm \ref{alg:HJ_computation_prob2} first solve the HJ PDEs (\eqref{eq:HJeq2_TI} or \eqref{eq:HJeq2_optCtrl_TI}) whose Hamiltonian has the maximum operation with 0.
Line 11 updates $V_2^\pm$ (or $V_2$) by choosing between $c$ and itself so that the HJ equation \eqref{eq:HJeq2_TI} (or \eqref{eq:HJeq2_optCtrl_TI}) holds.

\begin{algorithm}
\caption{Computing the solution $V_2^\pm$ or $V_2$ to the HJ equations for Problem \ref{prob:2} in the zero-sum game and optimal control settings. This algorithm deals with the four HJ equations: \eqref{eq:HJeq2}, \eqref{eq:HJeq2_TI}, \eqref{eq:HJeq2_opt}, and \eqref{eq:HJeq2_optCtrl_TI}.}
\begin{algorithmic}[1]
\State \textbf{Input:} {the temporal discretization: $\{t_0=0,...,t_K=T\}$, the spatial discretization: $\{(x_0,z_0),...,(x_N,z_N)\}$}
\State \textbf{Output:} {$V_2^\pm$ (or $V_2$)}
\State $V_2^\pm (\text{or }V_2) (T,x_i,z_i) \leftarrow \max\{c(T,x_i),g(T,x_i)-z_i\},\forall i$
\For {$k\in\{K-1,...,0\}$}
    \Case{solving the HJ equations \eqref{eq:HJeq2} or \eqref{eq:HJeq2_opt}}
        \State \Longunderstack[l]{$V_2^\pm (\text{or }V_2)(t_{k},x_i,z_i) \leftarrow V_2^\pm(\text{or }V_2) (t_{k+1},x_i,z_i)-$\\
        \quad\quad\quad $\Delta_k \hat{\bar{H}}^\pm (\text{or }\hat{\bar{H}}) (\phi_x^+,\phi_x^-,\phi_z^+,\phi_z^-), \forall i$}
        
        \State \Longunderstack[l]{$V_2^\pm (\text{or }V_2)(t_{k},x_i,z_i) \leftarrow \min\{g(t_k,x_i)-z_i,$\\
        \quad\quad\quad $V_2^\pm(\text{or }V_2) (t_{k},x_i,z_i)\}, \forall i$}
        
        \State \Longunderstack[l]{$V_2^\pm (\text{or }V_2)(t_{k},x_i,z_i) \leftarrow \max\{c(t_k, x_i),$\\
        \quad\quad\quad $V_2^\pm(\text{or }V_2) (t_{k},x_i,z_i)\}, \forall i$}
    \EndCase
    \Case{solving the HJ equations \eqref{eq:HJeq2_TI} or \eqref{eq:HJeq2_optCtrl_TI}}
        \State \Longunderstack[l]{$V_2^\pm(\text{or }V_2) (t_{k},x_i,z_i) \leftarrow V_2^\pm(\text{or }V_2) (t_{k+1},x_i,z_i)-$\\
        \quad\quad $\Delta_k \max\{0,\hat{\bar{H}}^\pm (\text{or }\hat{\bar{H}}) (\phi_x^+,\phi_x^-,\phi_z^+,\phi_z^-)\}, \forall i$}
        
        \State \Longunderstack[l]{$V_2^\pm (\text{or }V_2)(t_{k},x_i,z_i) \leftarrow \max\{c(t_k, x_i),$\\
        \quad\quad\quad $V_2^\pm(\text{or }V_2) (t_{k},x_i,z_i)\}, \forall i$}
    \EndCase
\EndFor
\end{algorithmic}
\label{alg:HJ_computation_prob2}
\end{algorithm}

%%%%%%%%%%%%%%%%%%%%%%%%%%%%%%%%%%%%%%%%%%%%%%%%%%%%%%%%%%%%%%%%%%%%%%%%%%%%%%%%
\section{Example} 
\label{sec:example}

\begin{figure}[t!]
\centering
\includegraphics[trim = 0mm 0mm 0mm 0mm, clip, width=0.35\textwidth]{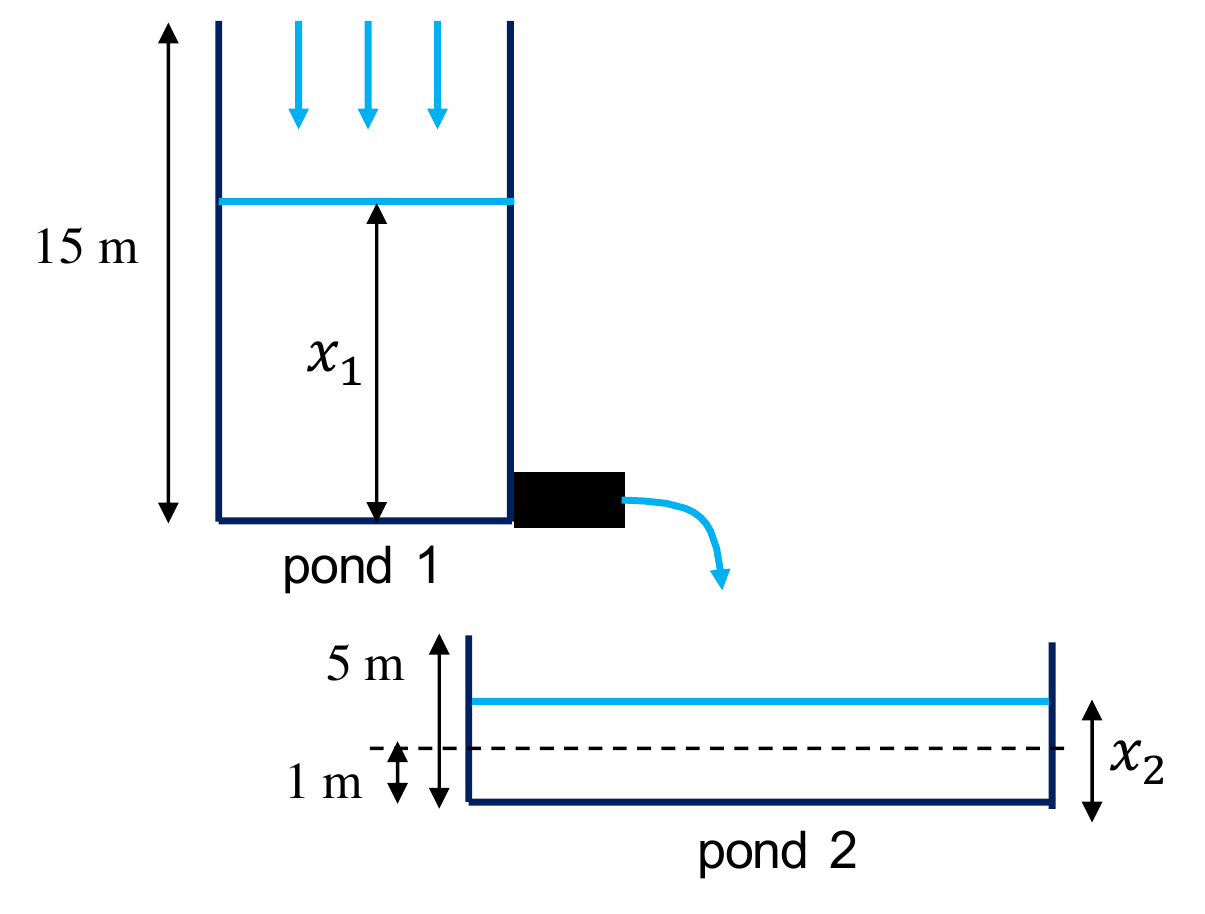} 
\caption{ Water system with two ponds.}
\label{fig:example}
\end{figure}

This section provides an example for Problem \ref{prob:1}; examples for Problem \ref{prob:2} can be found in \cite{donggun2020CDC}.
In this example, we solve a zero-sum game for two ponds, as shown in Figure \ref{fig:example}.
This example is motivated by the water system in \cite{chapman2018reachability}.

% Figure \ref{fig:example} shows a water system with two ponds.
Precipitation on pond 1 increases the water level of pond 1, and pond 1 (player A) wants to minimize the highest water level in the time horizon (1 s) by controlling amount of outflow to pond 2.
We assume that the water level increasing rate on pond 1 due to the precipitation is unknown but bounded by 0 and 10 $m/s$.
The precipitation is considered as player B.
These numbers and units can be easily changed to realistic problems.
We determine an optimal control signal and strategy even in the worst behavior of player B.

In the water system, we have two states: $x_1$ and $x_2$ represent the water level of pond 1 and 2.
The state trajectories are solving the following dynamics:
\begin{align}
\begin{split}
    &\dot{\mathrm{x}}_1(s) = \beta(s)-\sqrt{2g\mathrm{x}_1(s)} \alpha(s),\\
    &\dot{\mathrm{x}}_2(s) = 0.5\sqrt{2g\mathrm{x}_1(s)} \alpha(s)- 0.5\mathrm{x}_2(s),
\end{split}
\label{eq:example_dyn1}
\end{align}
where $\alpha(s)\in A=[0,1]$, $\beta(s)\in B=[0,10]$, and $g$ is the gravitational constant: 9.81 $m/s^2$.
In the dynamics for $\mathrm{x}_1$, the first term $\beta$ is by the precipitation (player B), and the second term $\sqrt{2g\mathrm{x}_1}\alpha$ is the water level decreasing rate by pond 1 (player A). 
The term $\sqrt{2g\mathrm{x}_1}$ is by Bernoulli's equation, and pond 1 controls the area of outflows ($\alpha$) between 0 and 1.
We set the bottom area of pond 2 is twice bigger than pond 1, thus the dynamics for $\mathrm{x}_2$ contains $0.5\sqrt{2g\mathrm{x}_1}\alpha$.
Also, we assume that pond 2's water is used for drinking water, which causes a decreasing rate $0.5\mathrm{x}_2$.

The dynamics \eqref{eq:example_dyn1} is not Lipschitz at $x_1=0$. To avoid this, we approximate $\sqrt{2gx_1}$ with a sinusoidal function $4.82\sin(1.17x_1)$ if $x_1$ is less than 1. 
This sinusoidal-approximate function has the same value and first derivative at $x_1=1$: $\sqrt{2g}\simeq 4.82\sin(1.17)$ and $\sqrt{g/2}\simeq 4.82*1.17*\cos(1.17)$.
The approximated (Lipschitz) dynamics are
\begin{align}
\begin{split}
    &\dot{\mathrm{x}}_1(s) = \beta(s)-\begin{cases}\sqrt{2g\mathrm{x}_1(s)} \alpha(s), &\mathrm{x}_1(s)\geq 1,\\ 4.82\sin(1.17\mathrm{x}_1(s))\alpha(s),&\mathrm{x}_1(s)<1,\end{cases} \\
    &\dot{\mathrm{x}}_2(s) = \begin{cases}0.5\sqrt{2g\mathrm{x}_1(s)} \alpha(s), &\mathrm{x}_1(s)\geq 1,\\ 2.41\sin(1.17\mathrm{x}_1(s))\alpha(s),&\mathrm{x}_1(s)<1,\end{cases}- 0.5\mathrm{x}_2(s),
\end{split}
\label{eq:example_dyn2}
\end{align}

We solve the two zero-sum games: the upper value function is 
\begin{align}
    &\quad\quad \vartheta_1^+ (0,x_1,x_2), =\max_{\delta\in\Delta(0)}\min_{\alpha\in\mathcal{A}(0)}\max_{\tau\in[0,1]} \mathrm{x}_1(\tau),\\
    &\text{subject to } \max\{|\mathrm{x}_1(s)-7.5|-7.5, |\mathrm{x}_2(s)-3|-2\}\leq 0,
    \label{eq:example_vartheta1+}
\end{align}
where $A=[0,1]$, $B=[0,10]$, $\mathcal{A}(0)=\{[0,1]\rightarrow A~|~\|\alpha\|_{L^\infty (0,1)}<\infty \}$, $\mathcal{B}(0)=\{[0,1]\rightarrow B ~|~ \|\beta\|_{L^\infty (0,1)}<\infty \}$, $\Delta(0)$ is a set of non-anticipative strategies for player B (pond 2) as in \eqref{eq:playerB_strategySet}, and $(\mathrm{x}_1,\mathrm{x}_2)$ solves \eqref{eq:example_dyn2} for $(\alpha,\delta[\alpha])$;
and the lower value function is
\begin{align}
    &\quad\quad \vartheta_1^- (0,x_1,x_2), =\min_{\gamma\in\Gamma(0)}\max_{\beta\in\mathcal{B}(0)}\max_{\tau\in[0,1]} \mathrm{x}_1(\tau),\\
    &\text{subject to } \max\{|\mathrm{x}_1(s)-7.5|-7.5, |\mathrm{x}_2(s)-3|-2\}\leq 0,
    \label{eq:example_vartheta1-}
\end{align}
where $\Gamma(0)$ is a set of non-anticipative strategies for player A (pond 1) as in \eqref{eq:playerA_strategySet}, and $(\mathrm{x}_1,\mathrm{x}_2)$ solves \eqref{eq:example_dyn2} for $(\gamma[\beta],\beta)$.
The state constraint implies that the water level of pond 1 has to be between 0 and 15 $m$ and the one of pond 2 is between 1 and 5 $m$.
In these games, pond 1 (player A) wants to minimize the worst water level of pond 1 in the time horizon while satisfying the state constraint for preventing flood in pond 1 and 2.

We will solve the HJ equation \eqref{eq:HJeq1} for $V_1^\pm$ corresponding to $\vartheta_1^\pm$ (\eqref{eq:example_vartheta1+} or \eqref{eq:example_vartheta1-}). 
We have the Hamiltonian
\begin{align}
% \begin{split}
    &\bar{H} ^+(t,x,z,p,q) = \max_{a\in A}\min_{b\in B}-p_1 b +0.5 p_2 x_2 \notag\\
    &\quad\quad + \begin{cases}(p_1 -0.5 p_2)\sqrt{2gx_1 } a &\text{if } x_1 \geq 1\\ (p_1 -0.5 p_2)4.82\sin(1.17 x_1 ) a &\text{if } x_1 <1\end{cases} \notag\\
    &\quad\quad = \begin{cases}-10p_1 &\text{if }p_1 \geq 0\\ 0& \text{if }p_1 <0\end{cases} + 0.5 p_2 x_2  \notag\\
    % & \begin{cases}(p_1 -0.5 p_2)\sqrt{2gx_1 } , & p_1 -0.5 p_2 \geq 0~\&~x_1 \geq 1,\\ (p_1 -0.5 p_2)4.82\sin(1.17 x_1 ),& p_1 -0.5 p_2  \geq 0~\&~x_1<1, \\ 0, & p_1-0.5p_2 <0\end{cases}
    % \\
    % & \begin{cases}(p_1 -0.5 p_2)\sqrt{2gx_1 } , & \begin{tabular}{l}$p_1 -0.5 p_2 \geq 0$\\ $x_1 \geq 1$\end{tabular}\\ (p_1 -0.5 p_2)4.82\sin(1.17 x_1 ),& p_1 -0.5 p_2  \geq 0~\&~x_1<1, \\ 0, & p_1-0.5p_2 <0\end{cases}
    % \\
    & \quad\quad +\begin{cases}(p_1 -0.5 p_2)\sqrt{2gx_1 } & \text{if }\begin{tabular}{l} $p_1 -0.5 p_2 \geq 0$\\ $x_1 \geq 1$\end{tabular}\notag\\ (p_1 -0.5 p_2)4.82\sin(1.17 x_1 ) &\text{if } \begin{tabular}{l} $p_1 -0.5 p_2  \geq 0$ \\ $x_1<1$\end{tabular} \\ 0 & \text{if } p_1-0.5p_2 <0\end{cases} \notag\\
    & = \bar{H} ^- (t,x,z,p,q)
% \end{split}
\label{eq:example_Hamiltonians}
\end{align}
where $x=(x_1,x_2)\in\R^2$ and $p=(p_1,p_2)\in\R^2$.
\eqref{eq:example_Hamiltonians} implies 
\begin{align}
    V_1^+ \equiv V_1^- \equiv V_1^\pm \text{ and } \vartheta_1^+\equiv \vartheta_1^- \equiv \vartheta_1^\pm.
\end{align}
We use $V_1^\pm$ to denote the same value functions $V_1^+$ and $V_1^-$.

% \begin{figure}[t!]
% \centering
% \includegraphics[trim = 0mm 0mm 0mm 0mm, clip, width=0.40\textwidth]{figures/ex1_ZeroLevelSet_0.5.pdf} 
% % \includegraphics[trim = 0mm 0mm 0mm 0mm, clip, width=0.40\textwidth]{figures/ex1_ZeroLevelSet_from_jpg.pdf} 
% \caption{ The zero-level set of $V_1^\pm$ is shown in this figure. The value of $V_1^\pm$ inside of the curvature is negative, but the value is positive outside.
% }
% \label{fig:example_zeroLevelSet}
% \end{figure}

\begin{figure}[t!]
\centering
\begin{tabular}{c}
     \includegraphics[trim = 0mm 0mm 0mm 0mm, clip, width=0.40\textwidth]{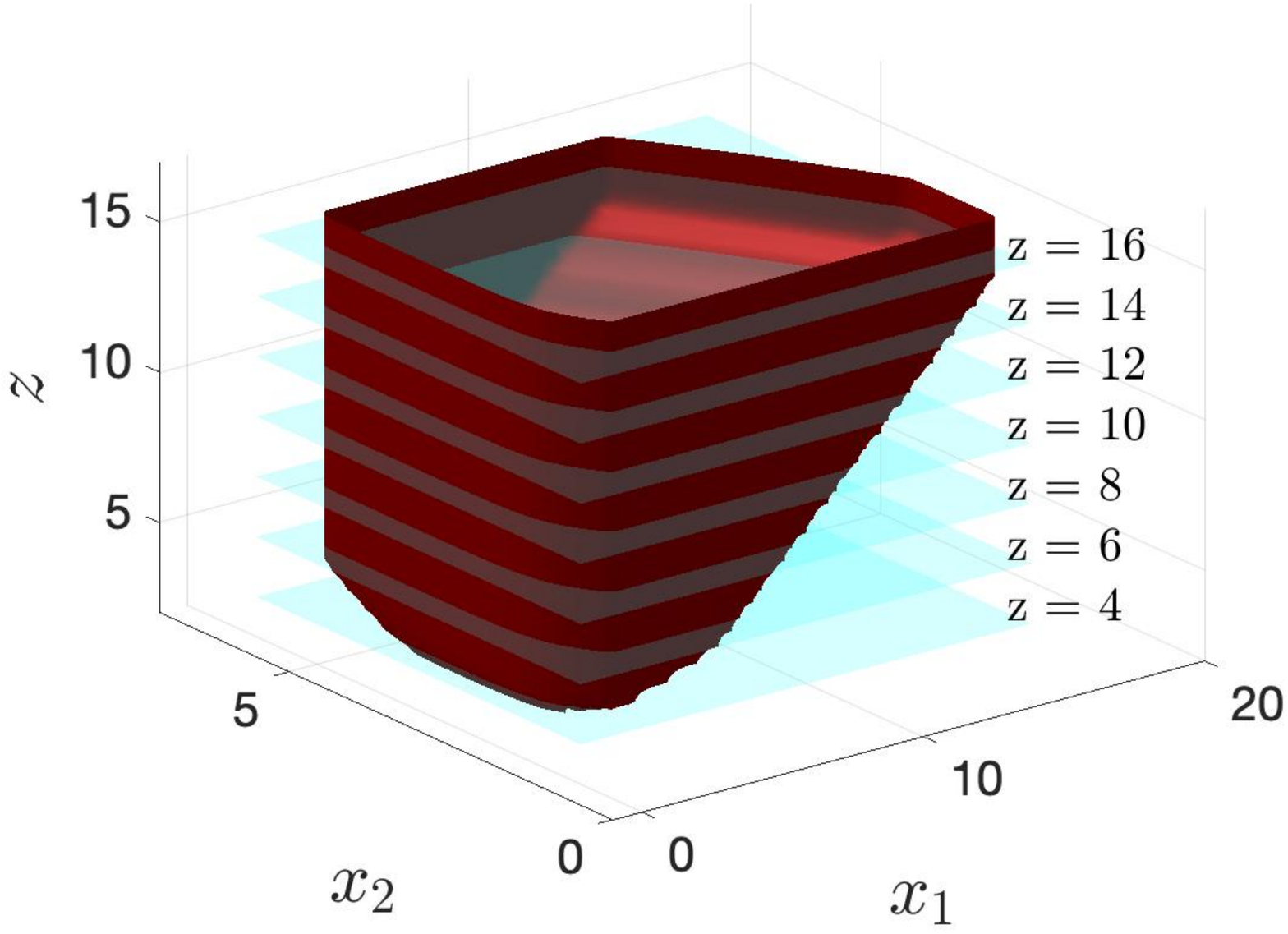}
    %  {figures/ex1_LevelSet_3D_from_jpg.pdf}
    %  {figures/ex1_LevelSet_3D.pdf} 
     \\(a)\\
     \includegraphics[trim = 0mm 0mm 0mm 0mm, clip, width=0.40\textwidth]{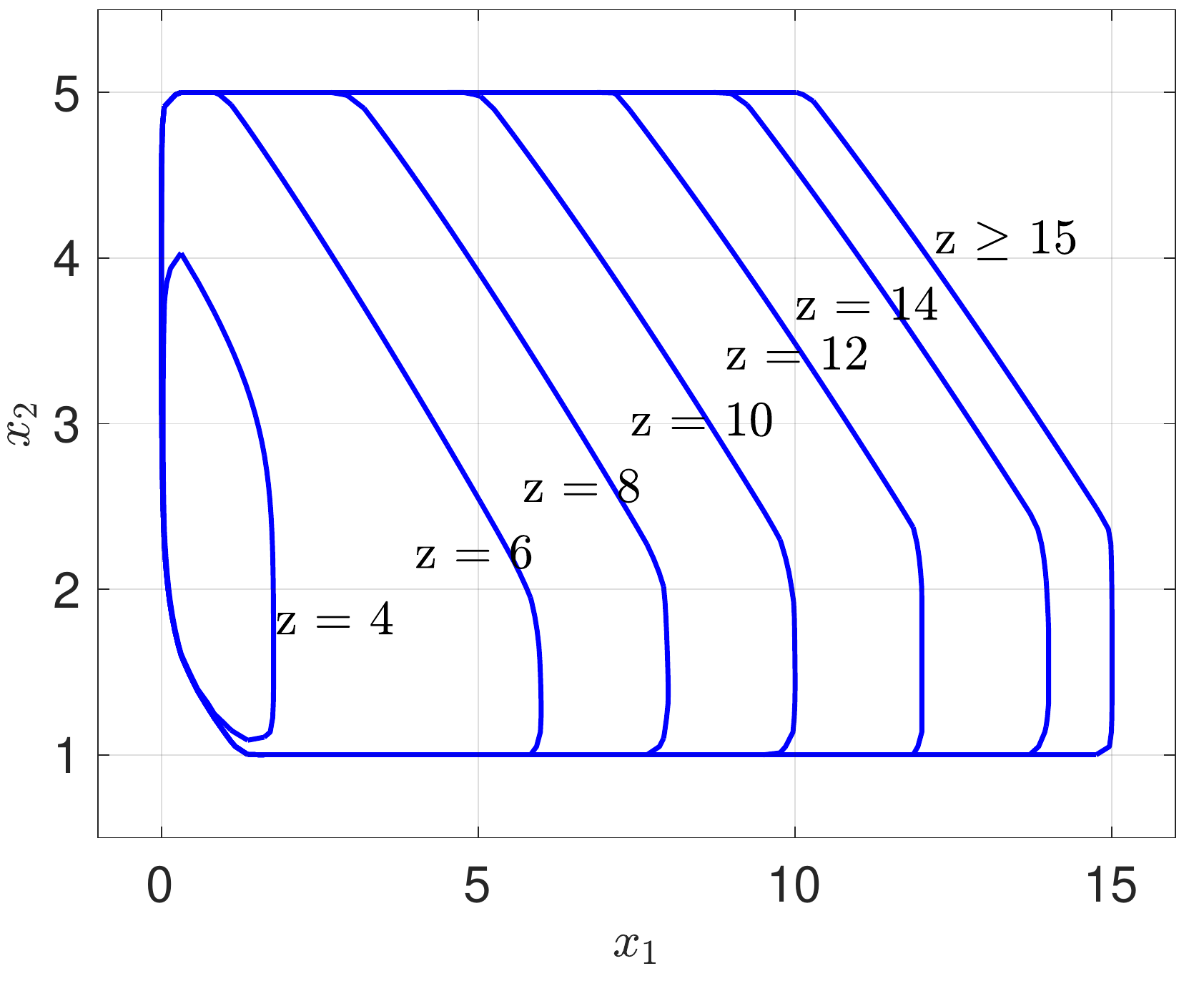}
    %  {figures/ex1_LevelSet_2D.pdf}  
     \\
    (b) 
\end{tabular}
\caption{ (a) The zero-level set of $V_1^\pm$ is shown in this figure. The value of $V_1^\pm$ inside of the curvature is negative, but the value is positive outside. The blue planes are $z$-level planes of 4, 6, 8, 10, 12, 14, and 16. (b) The $z$-level sets are shown in $(x_1,x_2)$-space. The $z$-level sets for $z\geq 15$ are the same (the outer curvature).
The $z$-level sets also show $\vartheta_1^\pm$ by Lemma \ref{lemma:Equiv_twoCost1}. 
For example, 
% for all $(x_1,x_2)$ on the $z$-level set of 4, $\vartheta_1^\pm$ is 4. Similarly, 
for $(1.60,2.85)$ on the $z$-level set of 4, $\vartheta_1^\pm$ is 4. 
On the other hand, consider $(0.11,2)$ on multiple $z$-level sets from 4.22 to any greather levels,
% This point is on multiple $z$-level sets from 4.22 to any greater levels, 
for which $\vartheta_1^\pm$ is the minimum $z$-level that contains the point: $\vartheta_1^\pm(0.05,2)=4.22$.
}
\label{fig:example_zLevelSet}
\end{figure}

The red curvature in Figure \ref{fig:example_zLevelSet} (a) shows the zero-level set of $V_1^\pm(0,x_1,x_2)$, numerically computed by Algorithm \ref{alg:HJ_computation_prob1}.
This algorithm is programmed by utilizing the \texttt{level set toolbox} \cite{mitchell2005toolbox} and the \texttt{helperOC toolbox} \cite{helperOC} in \texttt{Matlab}, and this simulation is carried out on a laptop with a 2.8 GHz Quad-Core i7 CPU and 16 GB RAM. 
Each of $x_1$, $x_2$, and $z$ axis has 81 discretization points, and the time interval $[0,1]$ is discretized with 201 points.
The computation time for $V_1^\pm$ is 237 $s$.
In Figure \ref{fig:example_zLevelSet} (a), the value of $V_1^\pm$ inside of the red curvature is negative, on the other hand, the value outside of the curvature is positive.

This example is time-invariant, so both HJ equations in \eqref{eq:HJeq1} and \eqref{eq:HJeq1_TI} can be utilized.
In this example, we solve the HJ equation \eqref{eq:HJeq1}.
% Note that the numerical solution to the time-invariant HJ equation \eqref{eq:HJeq1_TI} has very little difference between the solution to the HJ equation \eqref{eq:HJeq1} although we do not provide a graphical illustration.

Lemma \ref{lemma:Equiv_twoCost1} describes how to compute $\vartheta_1^\pm$ from the zero-level set of $V_1^\pm$, which is illustrated in Figure \ref{fig:example_zLevelSet}.
%the zero-level set of $V_1^\pm$ at each $z$-level 4 to 16.
Figure \ref{fig:example_zLevelSet} (a) illustrates the intersection of the zero-level set of $V_1^\pm$ and each $z$-level plane, and Figure \ref{fig:example_zLevelSet} (b) shows these intersections in the state space, $(x_1,x_2)$: the $z$-level sets on the zero-level set of $V_1^\pm$.
As illustrated in Figure \ref{fig:example_zLevelSet} (b), the lower $z$-level is achieved in the smaller region in $(x_1,x_2)$.
In this example, as $z$-level is increasing, the inner area of the $z$-level set on the subzero-level set of $V_1^\pm$ is increasing and also converging at the $z$-level of 15, which is the outer curvature in Figure \ref{fig:example_zLevelSet} (b).
For $(x_1,x_2)$ outside of the outer curvature indicated with $z\geq 15$, there is no control signal or strategy for pond 1 (player A) to satisfy the state constraint, which implies that $\vartheta_1^\pm (0,x_1,x_2)$ is infinity.
On the other hand, for $(x_1,x_2)$ on a unique $z$-level set, the $z$-level is equal to $\vartheta_1^\pm$. 
For example, the $z$-level set of 6 is the only $z$-level set passing through $(2.6,2)$. 
In this case, $\vartheta_1^\pm (0,2.6,2)=6$.
On the other hand, for $(x_1,x_2)$ on multiple $z$-level sets, the minimum value of $z$-level is $\vartheta_1^\pm$.
For example, $(0.05,2)$ is on the $z$-level sets of any number greater than or equal to 4.5.
In this case, $\vartheta_1^\pm (0,0.05,2)$ is 4.5 since 4.5 is the minimum $z$-level that contains the point $(0.05,2)$.

\begin{figure}[t!]
\centering
\begin{tabular}{c}
     \includegraphics[trim = 0mm 0mm 0mm 0mm, clip, width=0.40\textwidth]{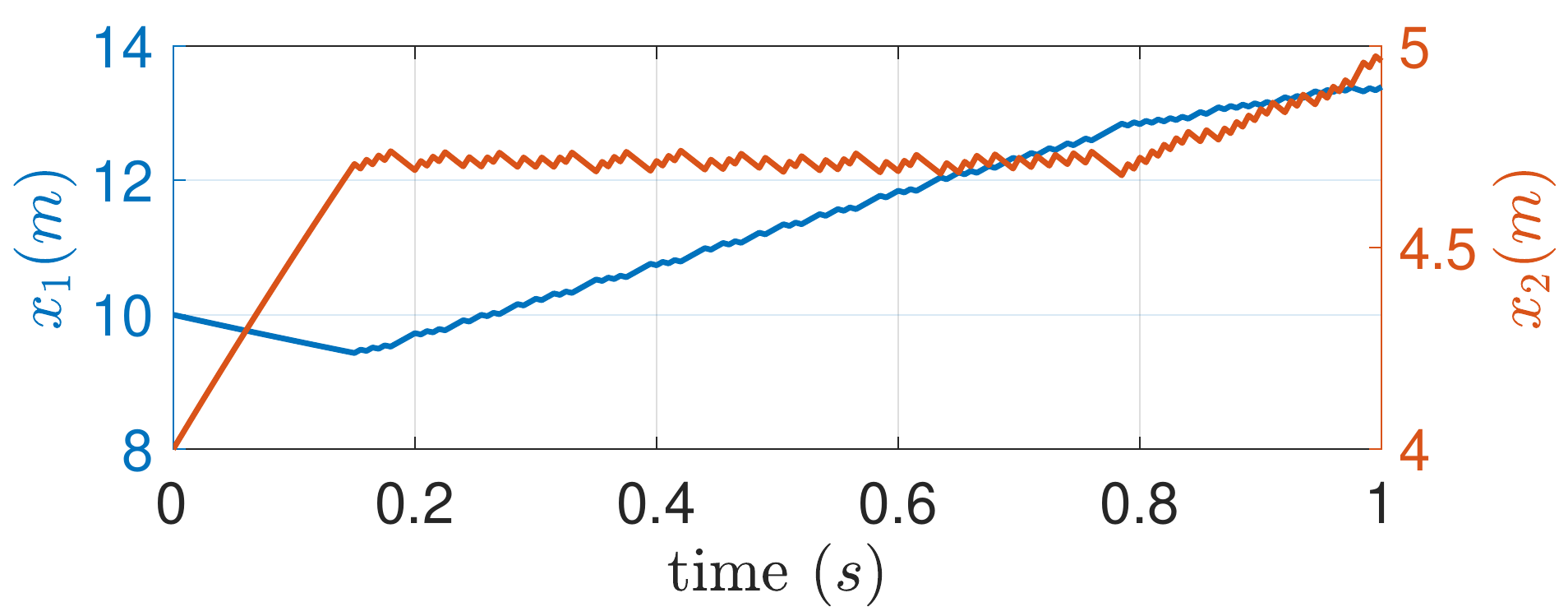}
     \\(a)\\
     \includegraphics[trim = 0mm 0mm 0mm 0mm, clip, width=0.40\textwidth]{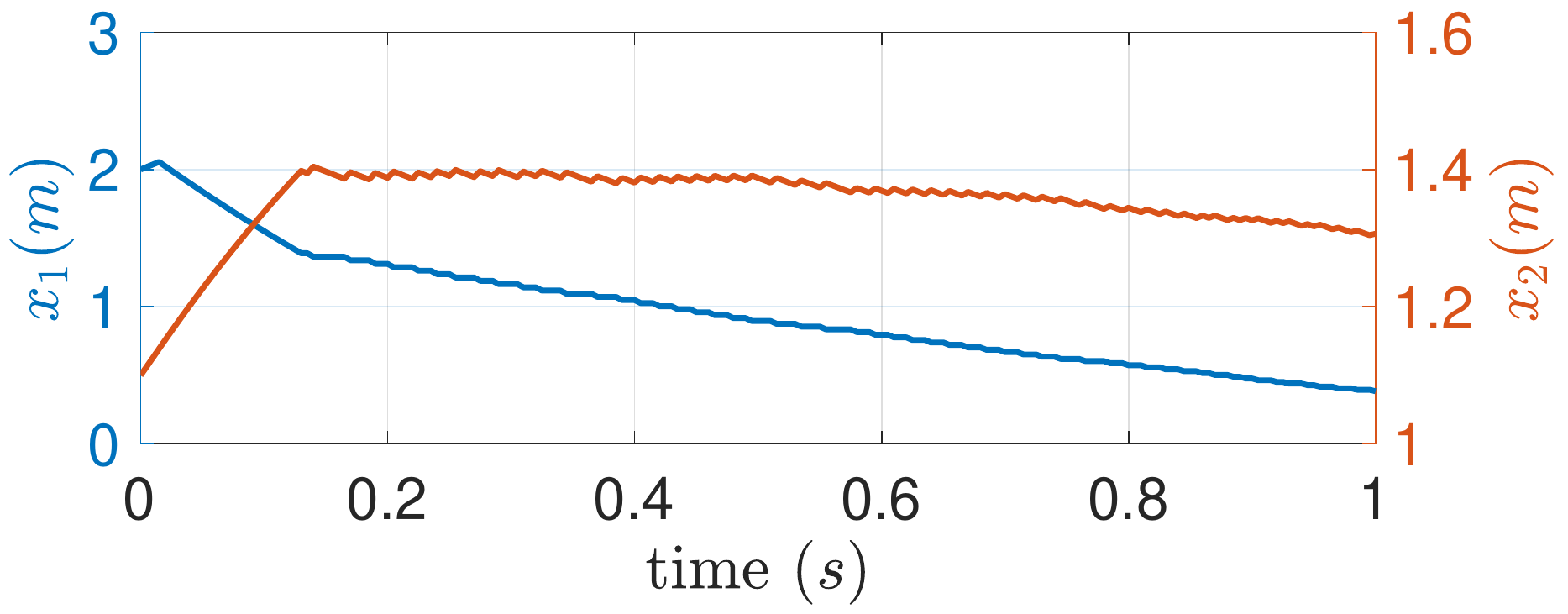}
    %  {figures/ex1_LevelSet_2D.pdf}  
     \\
    (b) 
\end{tabular}
\caption{ State trajectories by applying an optimal control signal and strategy for two players (pond 1 and the precipitation) where the initial states are (a) $(x_1,x_2)=(10,4)$ and (b) $(x_1,x_2)=(2,1.1)$.
}
\label{fig:optTraj}
\end{figure}

Using the value function $V_1^\pm$ and $\vartheta_1^\pm$, the method presented in Section \ref{sec:optCtrl_strategy} provides a state trajectory and an optimal control and strategy for the two players (pond 1 and the precipitation).
Among multiple solutions for optimal control and strategy presented in Remark \ref{remark:OptCtrl_Strat}, we choose
\begin{align}
\begin{split}
    &a_* \in \arg\max_{a\in A}\min_{b\in B}\tilde{H}_1^+ (a,b),\\
    &b_* \in \arg\min_{b\in B} \tilde{H}_1^+ (a_*,b),
\end{split}
\label{eq:bang-bang}
\end{align}
which satisfies the eight equations \eqref{eq:V1+_case1_a} to \eqref{eq:V1+_case2_b} since $\bar{H} ^+=\bar{H} ^-$ and $\tilde{H}_1^+ = \tilde{H}_1^-$, where $\tilde{H}_1^\pm$ is equal to $D_x V_1^\pm \cdot f +D_z V_1^\pm L$ as defined in \eqref{eq:tilde_Hamiltonian}.

Figure \ref{fig:optTraj} shows state trajectories for two different initial states: $(x_1,x_2)=(10,4)$ and $(2,1.1)$.
As shown in Figure \ref{fig:optTraj} (a), for the initial state $(10,4)$, $\mathrm{x}_2$ hits the boundary of the state constraint: $\mathrm{x}_2(1)=5$, and $\mathrm{x}_1$ is maximized at $t=1$.
Since the initial water levels of the two ponds are high, the precipitation (player B) tries to increase the water level of pond 1 for all time, but player A tries to balance the water levels of the two ponds.
On the other hand, for the initial state $(2,1.1)$, Figure \ref{fig:optTraj} (b) shows that $\mathrm{x}_2$ strictly satisfies the state constraint $[1,5]$. $\mathrm{x}_1$ is maximized at $t=0.015$ and increasing for the later time.
Since the initial water levels of the two ponds are low, the precipitation (player B) tries to violate the state constraint by not increasing the water level of pond 1.
However, player A tries to balance the two ponds' water level so that all ponds have more water than the minimum levels.
% It is shown, for the different initial states, that 
% the two terminal times that maximize the cost differs, which is not always the end of the time horizon.

As discussed in Section \ref{sec:NumericalAlgorithm}, there are some numerical issues in Algorithm \ref{alg:HJ_computation_prob1}.
% First, it is hard to get accurate gradients using level set methods \cite{Osher02}. 
% Since the optimal control or strategy are determined by the gradient, the state trajectories generated from these numerical control and strategy could have errors.
First, we observe that \eqref{eq:bang-bang} provides a bang-bang control, thus the state trajectories are not smooth as shown in Figure \ref{fig:optTraj}. This happens due to frequent sign change of the gradient along the time horizon.
Second, the numerical error on $V_1^\pm$ causes inaccurate $\vartheta_1^\pm$ by Lemma \ref{lemma:Equiv_twoCost1}, which could potentially cause unsafety even though the violation of the state constraint might be smaller for the smaller grid size.
In practice, we suggest having a safety margin to the state constraint: for example, use $c(s,\mathrm{x}(s))+\epsilon \leq 0$ for small $\epsilon >0$ instead of $c(s,\mathrm{x}(s))\leq0$.

%%%%%%%%%%%%%%%%%%%%%%%%%%%%%%%%%%%%%%%%%%%%%%%%%%%%%%%%%%%%%%%%%%%%%%%%%%%%%%%%
\section{Conclusion and Future Work} 
\label{sec:conclusion}

This paper presented four HJ equations for the two classes of  state-constrained zero-sum games where the terminal time is a variable to be determined and the stage cost is non-zero.
For each class of problems, two HJ equations have presented: one for time-varying version, and the other for the time-invariant version.
% Among these equations, four HJ equations are for the time-varying problems, and the other four HJ equations are for the time-invariant problems: the stage cost, the terminal cost, the dynamics, and the state constraint is time-invariant.
This paper also analyzed the optimal control and strategy for each player using the gradient of the viscosity solution to the HJ equations, and also presented a numerical algorithm to compute the viscosity solution.
As a practical example, a 2D water system demonstrates one of the presented HJ formulations.

% , and the computational complexity of the grid-based methods is discussed.

Although our HJ formulation can be generally utilized for the two classes of problems, the numerical computation of Algorithms \ref{alg:HJ_computation_prob1} and \ref{alg:HJ_computation_prob2} is intractable for high-dimensional system (higher than four or five). 
This is because the computational complexity is exponential in the dimension of the state.
In future, We aim to alleviate this complexity by deriving corresponding Lax and Hopf theory or by applying approximation techniques in reinforcement learning.

%%%%%%%%%%%%%%%%%%%%%%%%%%%%%%%%%%%%%%%%%%%%%%%%%%%%%%%%%%%%%%%%%%%%%%%%%%%%%%%%%%%%%%%%%%%%%%%%%%%%%%%%%%%%%%%%%%%%%%%%%%%%%%%%%%%%%%%%%%%%%%%%%%%%%%%%%%%%%%%%%%%%%%%%%%%%%%%%%%%%%%%%%%%%%%%%%%%%%%%%%%%%%%%%%%%%%%%%%%%%%%%%%%%%%%%%%%%%%%%%
\input{Appendix.tex}
%%%%%%%%%%%%%%%%%%%%%%%%%%%%%%%%%%%%%%%%%%%%%%%%%%%%%%%%%%%%%%%%%%%%%%%%%%%%%%%%%%%%%%%%%%%%%%%%%%%%%%%%%%%%%%%%%%%%%%%%%%%%%%%%%%%%%%%%%%%%%%%%%%%%%%%%%%%%%%%%%%%%%%%%%%%%%%%%%%%%%%%%%%%%%%%%%%%%%%%%%%%%%%%%%%%%%%%%%%%%%%%%%%%%%%%%%%%%%%%%
% \section*{Acknowledgements}
% The authors thank the name of people for discussions.

%%%%%%%%%%%%%%%%%%%%%%%%%%%%%%%%%%%%%%%%%%%%%%%%%%%%%%%%%%%%%%%%%%%%%%%%%%%%%%%%

\addtolength{\textheight}{0cm}   % This command serves to balance the column lengths

\bibliographystyle{IEEEtran}
\bibliography{IEEEexample}

\vspace*{-2\baselineskip}
\begin{IEEEbiography}
    [{\includegraphics[width=1in,height=1.25in,clip,keepaspectratio]{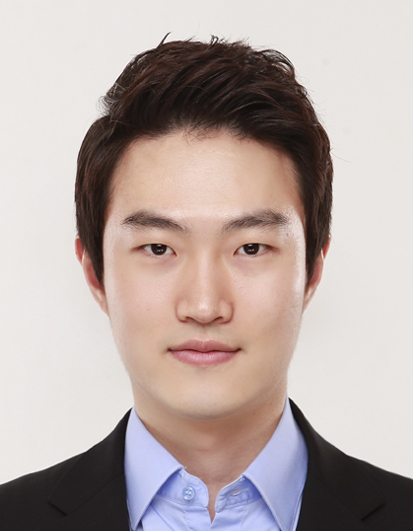}}]{Donggun Lee}
is a Ph.D. student in Mechanical Engineering at UC Berkeley. He received B.S. and M.S. degrees in Mechanical Engineering from Korea Advanced Institute of Science and Technology (KAIST), Daejeon, Korea, in 2009 and 2011, respectively. Donggun works in the area of control theory and robotics.
\end{IEEEbiography}
\vspace*{-2\baselineskip}
\begin{IEEEbiography}
    [{\includegraphics[width=1in,height=1.25in,clip,keepaspectratio]{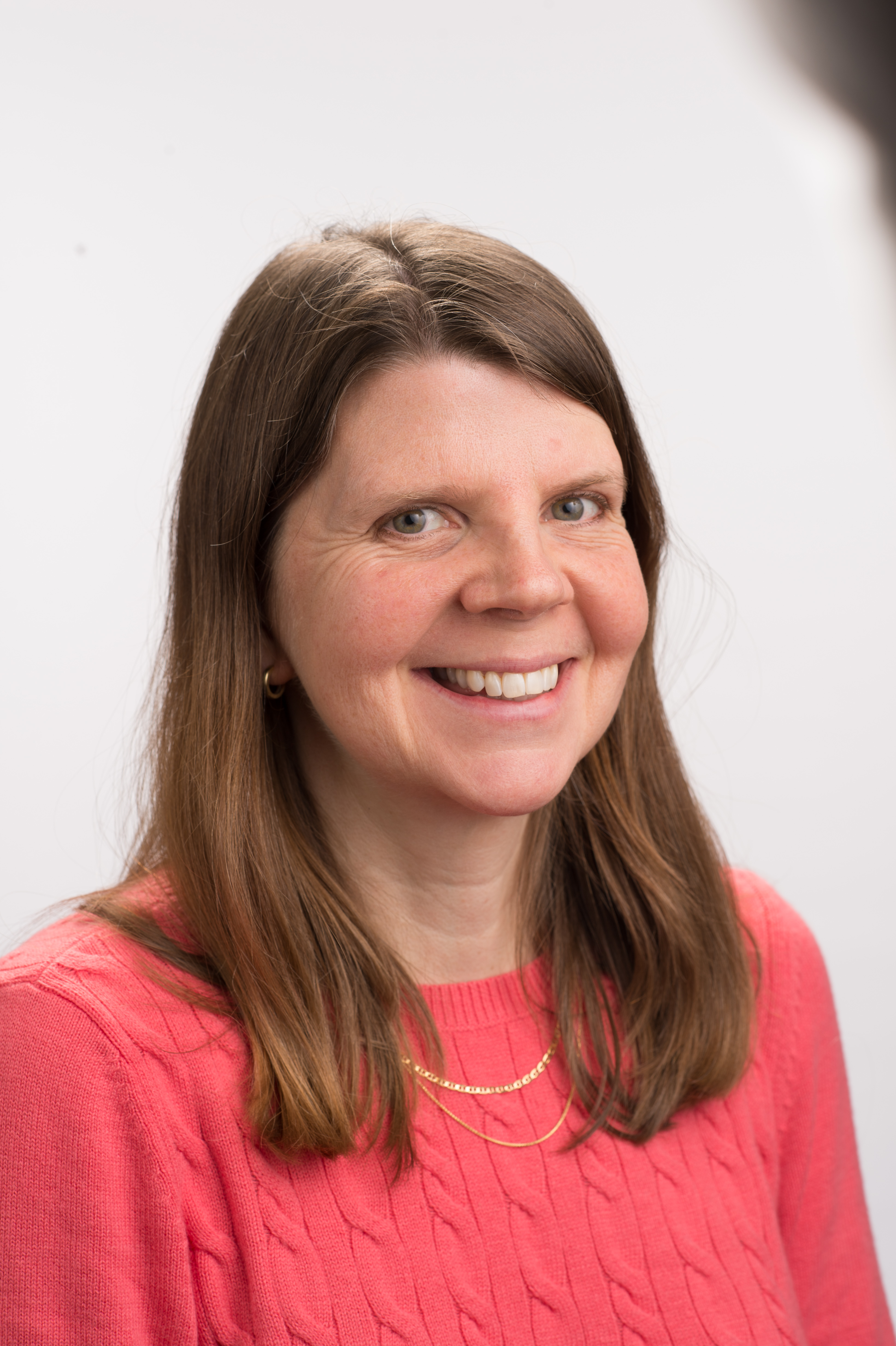}}]{Dr. Claire Tomlin}
is the Charles A. Desoer Professor of Engineering in EECS at Berkeley. She was an Assistant, Associate, and Full Professor in Aeronautics and Astronautics at Stanford from 1998 to 2007, and in 2005 joined Berkeley. Claire works in the area of control theory and hybrid systems, with applications to air traffic management, UAV systems, energy, robotics, and systems biology. She is a MacArthur Foundation Fellow (2006), an IEEE Fellow (2010), in 2017 she was awarded the IEEE Transportation Technologies Award, and in 2019 was elected to the National Academy of Engineering and the American Academy of Arts and Sciences.
\end{IEEEbiography}

\end{document}

%% file: Appendix.tex
\appendix
\subsection{Proof of Lemma \ref{lemma:Equiv_twoCost1}}
\label{appen:lemma_Equiv_twoCost1}
\textbf{Proof.} \\
(i) $\vartheta_1^+ (t,x) -z \leq0 \Rightarrow V_1^+ (t,x,z) \leq0$

$\vartheta_1^+ (t,x) -z\leq0$ implies that, for all $\delta\in\Delta(t)$, there exists $\alpha\in\mathcal{A}(t)$ such that 
\begin{align}
    \max_{\tau\in[t,T]}\int_t^\tau L(s,\mathrm{x}(s),\alpha(s),\delta[\alpha](s))ds+g(\tau,\mathrm{x}(\tau)) - z\leq \epsilon
\end{align}
and $c(s,\mathrm{x}(s))\leq0$ for $s\in[t,\tau]$ for any small $\epsilon>0$, where $\mathrm{x}$ solves \eqref{eq:dynamics} for $(\alpha,\delta[\alpha])$.
Thus, for all $\delta$, there exists $\alpha$ such that $J_1(t,x,z,\alpha,\delta) \leq \epsilon$, which concludes $V_1^+(t,x,z)\leq 0$.
Note that $J_1$ is defined in \eqref{eq:ctrl_costDef1}.

(ii) $V_1^+(t,x,z) \leq0  \Rightarrow \vartheta_1^+ (t,x) -z\leq 0$

Assumption \ref{assum:BigAssum} implies that, for any $\delta\in\Delta(t)$, there exists $\alpha\in\mathcal{A}$ such that $J_1(t,x,z,\alpha,\delta[\alpha])\leq V_1^+(t,x,z)$.
If $V_1^+(t,x,z)\leq0$, for any $\delta\in\Delta(t)$, there exists $\alpha$ such that $\max_{s\in[t,\tau]} c(s,\mathrm{x}(s)) \leq 0$ and $\int_t^\tau L(s,\mathrm{x}(s),\alpha(s),\delta[\alpha](s))ds + g(\tau,\mathrm{x}(s))-z \leq 0$ for all $\tau\in[t,T]$.
Thus, $\vartheta_1^+(t,x)-z\leq 0$.

(i) and (ii) concludes \eqref{eq:lemma_equiv_twoCost1}.

(iii) Let $\tilde\vartheta_1^+$ be the right hand term in \eqref{eq:vartheta1_upper_eq} subject to \eqref{eq:vartheta1_upper_eq_const}. 
Then, the following statement can be proved by analogous proofs in (i) and (ii).
\begin{align}
    \tilde\vartheta_1^+(t,x) = \min z &\text{ subject to } \sup_{\delta}\inf_{\alpha} \max\{ \max_{s\in[t,T]}c(s,\mathrm{x}(s)), \notag\\
    & \max_{\tau\in[t,T]}g(\tau,\mathrm{x}(\tau))-\mathrm{z}(\tau)\}\leq0.
\end{align}
By \eqref{eq:ctrl_costDef1_2} and \eqref{eq:lemma_equiv_twoCost1}, we conclude $\vartheta_1^+(t,x) = \tilde\vartheta_1^+(t,x)$.

(iv) The proof for $V_1^-$ and $\vartheta^-_1$ is similar to that for $V_1+$ and $\vartheta^+_1$.
\qed

\subsection{Proof of Lemma \ref{lemma:dynamicProgramming}}
\label{appen:lemma_dynamicProgramming}
\textbf{Proof.} Consider $(\mathrm{x},\mathrm{z})$ solving \eqref{eq:dynamics_comb} for any $(\alpha,\beta)$, and a small $h>0$. \eqref{eq:ctrl_costDef1_2} implies
\begin{align}
% \begin{split}
    J_1(t&,x,z,\alpha,\beta)=\max\Big\{ \max_{s\in[t,t+h]}c(s,\mathrm{x}(s)) ,\notag \\  &\max_{s\in[t,t+h]} g(s,\mathrm{x}(s))-\mathrm{z}(s),
    \max\big\{\max_{s\in[t+h,T]}c(s,\mathrm{x}(s)),\notag\\
    &\max_{s\in[t+h,T]} g(s,\mathrm{x}(s))-\mathrm{z}(s) \big\} \Big\}.
% \end{split}
\label{eq:lemma_dynprog_proof1}
\end{align}

(i) For all $\alpha\in\mathcal{A}(t)$ and $\delta\in\Delta(t)$, there exists $\alpha_1\in\mathcal{A}(t)$, $\delta_1\in\Delta(t)$, $\alpha_2\in\mathcal{A}(t+h)$, $\delta_2\in\Delta(t+h)$ such that
\begin{align}
    &\alpha(s) = \begin{cases} \alpha_1(s),&s\in[t,t+h],\\ \alpha_2(s),&s\in(t+h,T], \end{cases}\label{eq:lemma_dynprog_proof2}\\
% \end{align}
% \begin{align}
    &\delta[\alpha](s) = \begin{cases} \delta_1[\alpha](s),&s\in[t,t+h],\\ \delta_2[\alpha](s),&s\in(t+h,T]. \end{cases}
    \label{eq:lemma_dynprog_proof3}
\end{align}
Then, we have
\begin{align}
    &V_1^+(t,x,z) = \sup_{\substack{\delta_1\in\Delta(t)\\\delta_2\in\Delta(t+h)}} \inf_{\substack{\alpha_1\in\mathcal{A}(t)\\ \alpha_2\in\mathcal{A}(t+h) }} J_1(t,x,z,\alpha,\delta[\alpha]) \notag\\
    &= \sup_{\delta_1\in\Delta(t)}\inf_{\alpha_1\in\mathcal{A}(t)}\max\Big\{ \max_{s\in[t,t+h]}c(s,\mathrm{x}(s)) ,\notag \\  &\max_{s\in[t,t+h]} g(s,\mathrm{x}(s))-\mathrm{z}(s),
    \sup_{\delta_2\in\Delta(t+h)}\inf_{\alpha_2\in\mathcal{A}(t+h)}\max\notag\\
    &\big\{\max_{s\in[t+h,T]}c(s,\mathrm{x}(s)),\max_{s\in[t+h,T]} g(s,\mathrm{x}(s))-\mathrm{z}(s) \big\} \Big\}.
    \label{eq:lemma_dynprog_proof4}
\end{align}
The last equality is deduced by combining \eqref{eq:lemma_dynprog_proof1} and that the first two terms of $V_1^+$ ($\max_{s\in[t,t+h]}c(s,\mathrm{x}(s))$, $\max_{s\in[t,t+h]} g(s,\mathrm{x}(s))-\mathrm{z}(s)$) are independent of $(\alpha_2,\delta_2)$.
\eqref{eq:lemma_dynprog_proof4} concludes \eqref{eq:dynamicProgramming_V1_upper}.

(ii) The proof for \eqref{eq:dynamicProgramming_V1_lower} is similar to (i).
\qed

\subsection{Proof of Theorem \ref{thm:HJeq_prob1}}
\label{appen:thm_HJeq_prob1}
\textbf{Proof.}
\noindent
(i) At $t=T$, the definition of $V_1^\pm$ (\eqref{eq:game_V1_upper} and \eqref{eq:game_V1_lower}) implies \eqref{eq:V_terminalValue1}.

\noindent
(ii) For $U\in C^\infty ([0,T]\times\R^n\times\R)$ such that $V_1^+-U$ has a local maximum at $(t_0,x_0,z_0)\in(0,T)\times\R^n\times\R$ and $(V_1^+-U)(t_0,x_0,z_0)=0$, we will prove 
\begin{align}
\begin{split}
    \max\big\{&c(t_0,x_0)-U_0, g(t_0,x_0)-z_0-U_0, \\ 
    &U_{t0} - \bar{H}^+ (t_0,x_0,z_0,D_x U_0,D_z U_0)\big\}\geq0,
\end{split}
\label{eq:thm_HJeq_prob1_1}
\end{align}
where $U_0=U(t_0,x_0,z_0)$, $U_{t0}=U_t (t_0,x_0,z_0)$, $D_x U_{0}=D_x U (t_0,x_0,z_0)$, and $D_z U_0 = D_z U (t_0,x_0,z_0)$.

Suppose not. There exists $\theta>0$, $a_1\in A$ such that 
\begin{align}
    &c(t,x)-U_0<-\theta, \quad g(t,x)-z -U_0<-\theta , \\
\begin{split}
    &U_t(t,x,z) + D_x U(t,x,z)\cdot f(t,x,a_1,b) \\
    &\quad\quad\quad\quad\quad\quad - D_z U(t,x,z) L(t,x,a_1,b)   \leq -\theta
\end{split}
\end{align}
for all $b\in B$ and all points $(t,x,z)$ sufficiently close to $(t_0,x_0,z_0)$: $|t-t_0|+\|x-x_0\| + |z-z_0|<h_1$ for small enough $h_1>0$. 
Consider state trajectories $\mathrm{x}$ and $\mathrm{z}$ solving \eqref{eq:dynamics_comb} for $\alpha_1\equiv a_1$, $t=t_0$, $x=x_0$, $z=z_0$, and any $\beta\in\mathcal{B}(t_0)$.
By Assumption \ref{assum:BigAssum}, there exists a small $h$ such that $\|\mathrm{x}(s)-x_0\| + |\mathrm{z}(s)-z_0|<h_1-h$ ($s\in[t_0,t_0+h]$), then,
\begin{align}
% \begin{split}
    & c(s,\mathrm{x}(s))-U_0<-\theta,~ g(s,\mathrm{x}(s))-\mathrm{z}(s)-U_0<-\theta, \label{eq:thm_HJeq_prob1_2}\\
    % & g(s,\mathrm{x}(s))-\mathrm{z}(s)-U_0<-\frac{\theta}{2},
% \end{split}\\
% \end{align}
% and $\|\mathrm{x}(s)-x_0\|<\delta$ for all $s\in[t_0,t_0+h]$. Thus,
% \begin{align}
% \begin{split}
    &U_t(s,\mathrm{x}(s),\mathrm{z}(s)) \notag\\
    &+ D_x U(s,\mathrm{x}(s),\mathrm{z}(s))\cdot f(s,\mathrm{x}(s),a_1,\beta(s)) \notag\\
    &-D_z U(s,\mathrm{x}(s),\mathrm{z}(s)) L(s,\mathrm{x}(s),a_1,\beta(s))   \leq -\theta
% \end{split}
\label{eq:thm_HJeq_prob1_3}
\end{align}
for all $s\in[t_0,t_0+h]$ and $\beta\in\mathcal{B}(t_0)$.

Since $V_1^+-U$ has a local maximum at $(t_0,x_0,z_0)$,
\begin{align}
    &V_1^+(t_0+h,\mathrm{x}(t_0+h),\mathrm{z}(t_0+h))-V_1^+(t_0,x_0,z_0) \notag\\
    \leq &  U(t_0+h,\mathrm{x}(t_0+h),\mathrm{z}(t_0+h))-U(t_0,x_0,z_0)\notag\\
    =&  \int_{t_0}^{t_0+h} U_t(s,\mathrm{x}(s),\mathrm{z}(s)) \notag\\
    +& D_x U(s,\mathrm{x}(s),\mathrm{z}(s))\cdot f(s,\mathrm{x}(s),a_1,\delta[\alpha_1](s)) \notag\\
    -& D_z U(s,\mathrm{x}(s),\mathrm{z}(s)) L(s,\mathrm{x}(s),a_1,\delta[\alpha_1](s)) ds \leq -\theta h
    \label{eq:thm_HJeq_prob1_4}
\end{align}
for all $\delta\in\Delta(t_0)$, according to \eqref{eq:thm_HJeq_prob1_3}.
Lemma \ref{lemma:dynamicProgramming} implies
\begin{align}
    &V_1^+(t_0,x_0,z_0)\leq \sup_{\delta\in\Delta(t_0)} \max\big\{ \max_{s\in[t_0,t_0+h]} c(s,\mathrm{x}(s)),\notag\\
    &\max_{s\in[t,t+h]}g(\mathrm{x}(s))-\mathrm{z}(s),V_1^+(t_0+h,\mathrm{x}(t_0+h),\mathrm{z}(t_0+h))\big\}.
    \label{eq:thm_HJeq_prob1_5}
\end{align}
By subtracting $U_0$ on the both sides in \eqref{eq:thm_HJeq_prob1_5} and then applying \eqref{eq:thm_HJeq_prob1_2} and \eqref{eq:thm_HJeq_prob1_4}, we have
\begin{align}
    0\leq \max\{-\theta,-\theta,-\theta h\} <0,
\end{align}
which is contradiction. Thus, \eqref{eq:thm_HJeq_prob1_1} is proved.
\\
\noindent
(iii) For $U\in C^\infty ([0,T]\times\R^n\times\R)$ such that $V_1^+-U$ has a local minimum at $(t_0,x_0,z_0)\in(0,T)\times\R^n\times\R$ and $(V_1^+ -U)(t_0,x_0,z_0)=0$, we will prove 
\begin{align}
\begin{split}
    \max\big\{&c(t_0,x_0)-U_0, g(t_0,x_0)-z_0-U_0, \\ 
    &U_{t0} - \bar{H}^+ (t_0,x_0,z_0,D_x U_0,D_z U_0)\big\}\leq0,
\end{split}
\label{eq:thm_HJeq_prob1_6}
\end{align}

Since $J_1(t_0,x_0,z_0,\alpha)$ \eqref{eq:ctrl_costDef1} is greater than the value at $\tau=t_0$,
\begin{align}
    J_1(t_0,x_0,z_0,\alpha,\delta[\alpha]) \geq \max\big\{c(x_0,x_0),g(t_0,x_0)-z_0\big\},
\end{align}
for all $\alpha\in\mathcal{A}(t_0),\delta\in\Delta(t_0)$. By subtracting $U_0$ on the both sides, and taking the supremum over $\delta$ and the infimum over $\alpha$, sequentially, on the both side, we have
\begin{align}
    0 \geq \max\big\{c(x_0,x_0)-U_0,g(t_0,x_0)-z_0-U_0\big\}.
    \label{eq:thm_HJeq_prob1_7}
\end{align}

The rest of the proof is to show
\begin{align}
    U_{t0}-\bar{H}^+(t_0,x_0,z_0,D_xU_0,D_zU_0) \leq 0.
    \label{eq:thm_HJeq_prob1_10}
\end{align}
Suppose not. For some $\theta>0$, 
\begin{align}
\begin{split}
    U_t(t,x,z)  +&\max_{b\in B}D_x U(t,x,z)\cdot f(t,x,a,b) \\
    &-D_z U(t,x,z) L(t,x,a,b)  \geq \theta
\end{split}
\end{align}
for all $a\in A$ and all points $(t,x,z)$ sufficiently close to $(t_0,x_0,z_0)$: $|t-t_0|+\|x-x_0\| + |z-z_0|<h_1$ for small enough $h_1>0$.
Consider state trajectories $\mathrm{x}_1$ and $\mathrm{z}_1$ solving \eqref{eq:dynamics_comb} for any $\alpha\in\mathcal{A}(t_0)$, 
$\beta=\delta_1[\alpha]$, where
\begin{align}
    &\delta_1[\alpha](s)\in\arg\max_{b\in B}D_x U(s,\mathrm{x}_1(s),\mathrm{z}_1(s))\cdot f(s,\mathrm{x}_1(s),\alpha(s),b)\notag \\
    &\quad\quad\quad -D_z U(s,\mathrm{x}_1(s),\mathrm{z}_1(s)) L(s,\mathrm{x}_1(s),\alpha(s),b),
    \label{eq:thm_HJeq_prob1_11}
\end{align}
$t=t_0$, $x=x_0$, and $z=z_0$.
Since there exists a small $h>0$ such that $\|\mathrm{x}_1(s)-x_0\| + |\mathrm{z}_1(s)-z_0|<h_1-h$ ($s\in[t_0,t_0+h]$),
\begin{align}
\begin{split}
    &U_t(s,\mathrm{x}_1(s),\mathrm{z}_1(s)) \\
    +& D_x U(s,\mathrm{x}_1(s),\mathrm{z}_1(s))\cdot f(s,\mathrm{x}_1(s),\alpha(s),\delta_1[\alpha](s)) \\
    -&D_z U(s,\mathrm{x}_1(s),\mathrm{z}_1(s)) L(s,\mathrm{x}_1(s),\alpha(s),\delta_1[\alpha](s))   \geq \theta
\end{split}
\label{eq:thm_HJeq_prob1_8}
\end{align}
for all $s\in[t_0,t_0+h]$. 
By integrating \eqref{eq:thm_HJeq_prob1_8} over $s\in[t_0,t_0+h]$, we have
\begin{align}
    U(t_0+h,\mathrm{x}_1(t_0+h),\mathrm{z}_1(t_0+h)) - U(t_0,x,z) \geq \theta h.
    \label{eq:thm_HJeq_prob1_80}
\end{align}
Since \eqref{eq:thm_HJeq_prob1_80} holds for all $\alpha\in\mathcal{A}(t_0)$ and $\delta\in\Delta(t_0)$, 
\begin{align}
\begin{split}
    \sup_{\delta\in \Delta(t_0)}\inf_{\alpha\in\mathcal{A}(t_0)} U(t_0+h,\mathrm{x}(t_0+h),\mathrm{z}(t_0+h)) \\
    \quad\quad - U(t_0,x,z)   \geq \theta h,
\end{split}
\label{eq:thm_HJeq_prob1_21}
\end{align}
% \begin{align}
% \begin{split}
%     &\sup_{\delta\in \Delta(t_0)}\inf_{\alpha\in\mathcal{A}(t_0)} U_t(s,\mathrm{x}(s),\mathrm{z}(s)) \\
%     +& D_x U(s,\mathrm{x}(s),\mathrm{z}(s))\cdot f(s,\mathrm{x}(s),\alpha(s),\delta[\alpha](s)) \\
%     -&D_z U(s,\mathrm{x}(s),\mathrm{z}(s)) L(s,\mathrm{x}(s),\alpha(s),\delta[\alpha](s))   \geq \theta,
% \end{split}
% \label{eq:thm_HJeq_prob1_21}
% \end{align}
% for $s\in[t_0,t_0+h]$ 
where $\mathrm{x},\mathrm{z}$ solve \eqref{eq:dynamics_comb} for $(\alpha,\delta,t_0,x_0,z_0)$.

Since $V_1^+-U$ has a local minimum at $(t_0,x_0,z_0)$,
\begin{align}
    &\sup_{\delta\in \Delta(t_0)}\inf_{\alpha\in\mathcal{A}(t_0)}\begin{tabular}{l}$V_1^+(t_0+h,\mathrm{x}(t_0+h),\mathrm{z}(t_0+h))$\\
    $\quad\quad\quad\quad\quad\quad\quad\quad\quad -V_1^+(t_0,x_0,z_0)$ \end{tabular}\notag \\
    \geq &\sup_{\delta\in \Delta(t_0)}\inf_{\alpha\in\mathcal{A}(t_0)} \begin{tabular}{l} $U(t_0+h,\mathrm{x}(t_0+h),\mathrm{z}(t_0+h))$\\
    $\quad\quad\quad\quad\quad\quad\quad\quad\quad-U(t_0,x_0,z_0)$\end{tabular} \notag\\
    \geq & ~\theta h \label{eq:thm_HJeq_prob1_9}
    % =& \sup_{\delta\in \Delta(t_0)}\inf_{\alpha\in\mathcal{A}(t_0)}\int_{t_0}^{t_0+h} U_t(s,\mathrm{x}(s),\mathrm{z}(s)) \notag\\
    % & + D_x U(s,\mathrm{x}(s),\mathrm{z}(s))\cdot f(s,\mathrm{x}(s),\alpha(s),\delta[\alpha](s)) \label{eq:thm_HJeq_prob1_9}\\
    % & -D_z U(s,\mathrm{x}(s),\mathrm{z}(s)) L(s,\mathrm{x}(s),\alpha(s),\delta[\alpha](s)) ds \geq \theta h \notag
\end{align}
according to \eqref{eq:thm_HJeq_prob1_21}. However, Lemma \ref{lemma:dynamicProgramming} implies
\begin{align}
% \begin{split}
    \sup_{\delta\in\Delta(t_0)}\inf_{\alpha\in\mathcal{A}(t_0)}V_1^+ (t_0+h, \mathrm{x}&(t_0+h), \mathrm{z}(t_0+h)) \notag\\ &\leq V_1^+(t_0,x_0,z_0),
% \end{split}
\end{align}
which contradicts \eqref{eq:thm_HJeq_prob1_9}.

(iv) The proof for the viscosity solution $V_1^-$ is similar to (ii) and (iii) for $V_1^+$.
Also, the uniqueness follows from the uniqueness theorems for viscosity solutions, Theorem 4.2 in \cite{barron1989bellman}, and the extension of Theorem 1 in \cite{evans10}.
\qed

%%%%%%%%%%%%%%%%%%%%%%%%%%%%%%%%%%%%%%%%%%%%%%%%%%%%%%%%%%%%%%%%%%%%%%%%%%%%%%%%%%%%%%%%%%%%%%%%%%%%%%%%%%%%%%%%%%%%%%%%%%%%%%%%%%%%%%%%%%%%%%%%%%%%%%%%%%%%%%%%%%%%%
\subsection{Proof of Lemma \ref{lemma:equiv_V1}}
\label{appen:lemma_equiv_V1}

% \begin{align}
%     &V_1^+(t,x,z) =    \sup_{\substack{\delta\in\Delta(t),\\\nu_A\in\Delta_d(t)}}\inf_{\alpha\in\mathcal{A}(t)} \tilde{J}(t,x,z,\alpha,\delta[\alpha],\nu_A[\alpha]),\\
%     % \label{eq:game_V1_upper_finiteHorizon}\\
% % \end{align}
% % and 
% % \begin{align}
%     &V_1^-(t,x,z) =    \inf_{\tilde\gamma\in\tilde{\Gamma}(t)}\sup_{\substack{\beta\in\mathcal{B}(t),\\\mu\in\mathcal{M}(t)}} \tilde{J}(t,x,z,\tilde\gamma[\beta,\mu],\beta,\mu).
%     % \label{eq:game_V1_lower_finiteHorizon}
% \end{align}

\textbf{Proof.} Set $\tilde V_1^+$ and $\tilde V_1^-$ be the right hand terms in \eqref{eq:game_V1_upper_finiteHorizon} and \eqref{eq:game_V1_lower_finiteHorizon}, respectively.
$V_1^+$ are $V_1^-$ are defined in \eqref{eq:game_V1_upper} and \eqref{eq:game_V1_lower}, respectively.
% Set $\tilde{\vartheta}_1^+$ is \eqref{eq:vartheta1_upper_opt_cost} subject to \eqref{eq:vartheta1_upper_opt_const}, and $\vartheta_1^+$ is defined in \eqref{eq:vartheta1_upper_eq} subject to \eqref{eq:vartheta1_upper_eq_const}.

(i) In this proof, we utilize the following properties in \cite{mitchell2005time,lee2020hopf}, presented as below.

Define a pseudo-time operator $\sigma_\mu:[t,T]\rightarrow [t,T]$ for a given $\mu\in\mathcal{M}(t)$ (defined in \eqref{eq:ctrlSet_freezing}) and the corresponding inverse operator:
\begin{align}
    &\sigma_\mu(s  ) = \int_t^s \mu(\tau) d\tau + t;\\
% \end{align}
% and the corresponding inverse operator: 
% \begin{align}
    &\sigma^{-1}_\mu(s) \coloneqq \min \tau \text{ subject to } \sigma_\mu (\tau) = s.
\end{align}
Then, 
\begin{align}
    & \sigma_\mu\big( \sigma^{-1}_\mu(s) \big) = s,\quad s\in[t,\sigma_\mu(T)],\label{eq:lemma5_proof_00}\\
    & \sigma^{-1}_\mu\big( \sigma_\mu(s) \big) = s, \quad s\in\textrm{Range}(\sigma_\mu^{-1}),\label{eq:lemma5_proof_01}
\end{align}
where $\textrm{Range}(\sigma_\mu^{-1})\coloneqq \{ \sigma^{-1}_\mu(s)~|~s\in[t,\sigma_\mu(T)] \}$.

Consider two state trajectories: $(\mathrm{x},\mathrm{z})$ solving \eqref{eq:dynamics_comb} for $(\tilde\alpha(\sigma_\mu^{-1}(\cdot)), \tilde\beta(\sigma_\mu^{-1}(\cdot)))$ for $s\in[t,\sigma_\mu(T)]$; $(\tilde{\mathrm{x}},\tilde{\mathrm{z}})$ solving \eqref{eq:dynamics_comb_TI} for $(\tilde\alpha,\tilde\beta,\mu)$, and $\mathrm{x}(t)=\tilde{\mathrm{x}}(t) = x$.
Then,
\begin{align}
    & \mathrm{x} \big( \sigma_\mu(s) \big) = \tilde{\mathrm{x}}( s), \quad s\in[t,T],
    \label{eq:lemma5_proof0_1}\\
    % & \int_t^{\sigma_\mu(T) } L\big( \mathrm{x}(s), \tilde{\alpha }( \sigma^{-1}_\mu(s) ), \tilde{\beta }( \sigma^{-1}_\mu(s) )  \big)ds + g\big(\mathrm{x}(\sigma_\mu(T))\big) \notag\\
    % =& \int_t^T L\big( \tilde{\mathrm{x}}(s), \tilde{\alpha}(s), \tilde{\beta}(s) \big) \mu (s) ds + g\big(\tilde{\mathrm{x}}(T)\big).
    &  g\big(\mathrm{x}(\sigma_\mu(T))\big) - \mathrm{z}(\sigma_\mu(T)) =  g\big(\tilde{\mathrm{x}}(T)\big)-\tilde{\mathrm{z}}(T).
    \label{eq:lemma5_proof0_2}
\end{align}
\eqref{eq:lemma5_proof0_1} is according to Lemma 4 in \cite{Mitchell03}, and \eqref{eq:lemma5_proof0_2} is derived by combining two lemmas (Lemma 4 and 6) in \cite{Mitchell03}.

(ii) $\tilde{V}_1^+(t,x,z) \geq V_1^+(t,x,z)$

For small $\epsilon>0$, there exists $\delta_1\in\Delta(t)$ such that
\begin{align}
\begin{split}
    V_1^+ &(t,x,z)-\epsilon  \leq  \inf_{\alpha}\max_{\tau\in[t,T]}\max\big\{ \max_{s\in[t,\tau]}c(\mathrm{x}_1(s)),\\
    &\quad\quad\quad\quad\quad\quad\quad\quad\quad\quad\quad g\big(\mathrm{x}_1(\tau)\big)-\mathrm{z}_1(\tau)\big\},
\end{split}
\label{eq:lemma5_proof1}
\end{align}
where $(\mathrm{x}_1,\mathrm{z}_1)$ solves \eqref{eq:dynamics_comb} for $(\alpha,\delta_1[\alpha])$.
Denote $\tau_*(\alpha)$ is the maximizer of the right hand term in \eqref{eq:lemma5_proof1} for each $\alpha\in\mathcal{A}(t)$:
\begin{align}
\begin{split}
    &\tau_*(\alpha) \coloneqq \arg\max_{\tau\in[t,T]} \max\big\{ \max_{s\in[t,\tau]}c(\mathrm{x}_1(s)), g\big(\mathrm{x}_1(\tau)\big)-\mathrm{z}_1(\tau)\big\}.
\end{split}
\end{align}
Define a particular strategy $\nu_{A,1}\in\textrm{N}_A(t)$:
\begin{align}
    \nu_{A,1}[\alpha](s)\coloneqq \begin{cases} 1, & s\in[t,\tau_*(\alpha)], \\ 0,& s\in(\tau_*(\alpha),T]. \end{cases}
    \label{eq:lemma5_proof2}
\end{align}

Consider a state trajectory $(\tilde{\mathrm{x}}_1,\tilde{\mathrm{z}}_1)$ solving \eqref{eq:dynamics_comb_TI} for $(\alpha,\delta_1[\alpha],\nu_{A,1}[\alpha])$.
Then, we have 
\begin{align}
    &(\tilde{\mathrm{x}}_1,\tilde{\mathrm{z}}_1)(s) = \begin{cases} (\mathrm{x}_1,\mathrm{z}_1)(s), & s\in [t,\tau_*(\alpha)], \\ (\mathrm{x}_1,\mathrm{z}_1)(\tau_*(\alpha)), & s\in (\tau_*(\alpha),T],
    \end{cases} 
    \label{eq:lemma5_proof3}
\end{align}
% $$ for $s\in$, and $\mathrm{x}_1(\tau_*(\alpha)) = \mathrm{x}_2(T)$.
Since $\tilde{V}_1^+$ has the supremum over $(\delta,\nu_A)$-space operation,
\begin{align}
    &\tilde{V}_1^+(t,x,z) \geq \inf_{\alpha} \max\big\{ \max_{s\in[t,T]}c(\tilde{\mathrm{x}}_1(s)), g\big(\tilde{\mathrm{x}}_1(T)\big)- \tilde{\mathrm{z}}_1(T) \big\} \notag\\
    &= \inf_{\alpha} \max\big\{ \max_{s\in[t,\tau_*(\alpha)]}c( \mathrm{x}_1(s)),  g\big({\mathrm{x}}_1(\tau_*(\alpha))\big)-\mathrm{z}_1(\tau_*(\alpha)) \big\} \notag\\
    & \geq V_1^+ (t,x,z)-\epsilon. 
    \label{eq:lemma5_proof4}
\end{align}
The second equality is according to \eqref{eq:lemma5_proof3}, and the third inequality is by \eqref{eq:lemma5_proof1}.

%%%%%%%%%%%%%%%%%%%%%%%%%%%%%%%%%%%%%%%%%%
(iii) $V_1^+(t,x,z) \geq \tilde{V}_1^+(t,x,z)$

Define $\tilde{\mathfrak{A}}_\mu:\mathcal{A}(t)\rightarrow\mathcal{A}(t)$ and its psuedo inverse function ${\mathfrak{A}}_\mu:\mathcal{A}(t)\rightarrow\mathcal{A}(t)$:
\begin{align}
    & (\tilde{\mathfrak{A}}_\mu(\alpha))(s) \coloneqq \begin{cases} \alpha\big( \sigma_\mu (s) \big),
    & s\in\textrm{Range}(\sigma^{-1}_\mu), \\
    \text{any }a \in A, & s\notin\textrm{Range}(\sigma^{-1}_\mu), \end{cases}
    \label{eq:lemma5_proof3_1}\\
    &(\mathfrak{A}_\mu(\tilde{\alpha}))(s) \coloneqq \begin{cases} \tilde{\alpha}(\sigma^{-1}_\mu(s)), & s\in[t,\sigma_\mu(T)], \\\text{any }a\in A ,&s\in(\sigma_\mu(T),T],  \end{cases}
    \label{eq:lemma5_proof3_2}
\end{align}
Also, define $\tilde{\mathfrak{D}}_\mu:\Delta(t)\rightarrow\Delta(t)$ and its psuedo inverse function $\mathfrak{D}_\mu:\Delta(t)\rightarrow\Delta(t)$:
\begin{align}
    & (\tilde{\mathfrak{D}}_\mu(\delta))[\tilde{\alpha}](s) = \begin{cases} \delta[\mathfrak{A}_\mu(\tilde{\alpha}) ] ( \sigma_\mu(s) ) ,& s\in\textrm{Range}(\sigma^{-1}_\mu),\\ \text{any }b\in B, &s\notin\textrm{Range}(\sigma^{-1}_\mu), \end{cases} 
    \label{eq:lemma5_proof3_3}\\
    & (\mathfrak{D}_\mu(\tilde\delta))[\alpha](s) = \begin{cases} \tilde{\delta}[\tilde{\mathfrak{A}}_\mu(\alpha)](\sigma^{-1}_\mu(s)), & s\in[t,\sigma_\mu(T)], \\\text{any } b\in B, & s\in (\sigma_\mu(T),T]. \end{cases}
    \label{eq:lemma5_proof3_4}
\end{align}
These definitions satisfy the following properties:
\begin{align}
    \begin{split}
        &(\tilde{\mathfrak{A}}_\mu(\mathfrak{A}_\mu(\tilde\alpha)))(s)= \tilde{\alpha}(s),\\ &(\tilde{\mathfrak{D}}_\mu(\mathfrak{D}_\mu(\tilde\delta)))[\tilde\alpha](s) = \tilde\delta [\tilde\alpha](s), 
    \end{split}
    \quad \text{for } s\in\textrm{Range}(\sigma^{-1}_\mu)\\
    \begin{split}
        &(\mathfrak{A}_\mu(\tilde{\mathfrak{A}}_\mu(\alpha)))(s)= {\alpha}(s),\\ &(\mathfrak{D}_\mu(\tilde{\mathfrak{D}}_\mu(\delta)))[\alpha](s) = \delta [\alpha](s), 
    \end{split}
    \quad \text{for } s\in[t,\sigma_\mu(T)],\\
    % \\
    % & (\tilde{\mathfrak{D}}(\mathfrak{D}(\tilde\delta)))[\tilde\alpha](s) = (\mathfrak{D}(\tilde\delta))[\mathfrak{A}(\tilde\alpha)](\sigma(s;\beta_d)) \notag\\
    % & = \tilde\delta [\tilde{\mathfrak{A}}(\mathfrak{A}(\tilde\alpha))](s) = \tilde\delta [\tilde\alpha](s) 
% \end{align}
% Then,
% \begin{align}
    & \big\{ \alpha={\mathfrak{A}}_\mu(\tilde\alpha)~|~\tilde\alpha\in\mathcal{A}(t) \big\} = \mathcal{A}(t), \forall \mu\in\mathcal{M}(t)
    \label{eq:lemma5_proof3_5}\\
    & \big\{ \delta={\mathfrak{D}}_\mu(\tilde\delta)~|~\tilde\delta\in\Delta(t) \big\} = \Delta(t), \forall \mu\in\mathcal{M}(t).
    \label{eq:lemma5_proof3_6}
\end{align}

Consider $(\tilde{\mathrm{x}},\tilde{\mathrm{z}})$ solving \eqref{eq:dynamics_comb_TI} for $(\tilde\alpha,\tilde\delta[\tilde\alpha],\mu)$,  $(\mathrm{x},\mathrm{z})$ solving \eqref{eq:dynamics_comb} for $(\mathfrak{A}_\mu(\tilde\alpha),(\mathfrak{D}_\mu(\tilde\delta))[\mathfrak{A}_\mu(\tilde\alpha)])$, and $(\mathrm{x}_1,\mathrm{z}_1)$ solving \eqref{eq:dynamics_comb} for $(\alpha,\delta[\alpha])$.
Then, we have
\begin{align}
    &\sup_{\tilde\delta\in\Delta(t)} \inf_{\tilde\alpha\in\mathcal{A}(t)} \max\big\{ \max_{s\in[t,T]}c(\tilde{\mathrm{x}}(s)),g(\tilde{\mathrm{x}}(T))-\tilde{\mathrm{z}}(T)\big\} \notag\\
    =  &\sup_{\tilde\delta\in\Delta(t)} \inf_{\tilde\alpha\in\mathcal{A}(t)} \max\big\{ \max_{s\in[t,\sigma_\mu(T)]}c(\mathrm{x}(s)),\notag\\
     & \quad\quad\quad\quad\quad\quad\quad\quad g(\mathrm{x}(\sigma_\mu(T)))-{\mathrm{z}}(\sigma_\mu(T))\big\},\label{eq:lemma5_proof3_7}\\
    = & \sup_{\delta\in\Delta(t)} \inf_{\alpha\in\mathcal{A}(t)} \max\big\{ \max_{s\in[t,\sigma_\mu(T)]}c(\mathrm{x}_1(s)),\notag\\
    & \quad\quad\quad\quad\quad\quad\quad\quad g(\mathrm{x}_1(\sigma_\mu(T)))-{\mathrm{z}}_1(\sigma_\mu(T))\big\},
    \label{eq:lemma5_proof3_8}\\
    \leq & V_1^+(t,x,z). \label{eq:lemma5_proof3_9}
    % \sup_{\delta\in\Delta(t)} \inf_{\alpha\in\mathcal{A}(t)}\max_{\tau\in[t,T]} \max\big\{ \max_{s\in[t,\tau]}c(\mathrm{x}_1(s)),\notag\\ & \quad\quad\quad\quad\quad\quad\quad\quad g(\mathrm{x}_1(\tau))-{\mathrm{z}}_1(\tau)\big\}
\end{align}

\eqref{eq:lemma5_proof3_7} is by \eqref{eq:lemma5_proof0_1} and \eqref{eq:lemma5_proof0_2}, and \eqref{eq:lemma5_proof3_8} is according to \eqref{eq:lemma5_proof3_5} and \eqref{eq:lemma5_proof3_6}.
Since the above inequality holds for all $\mu$, we substitute $\nu_A[\alpha]$ for $\mu$ and take the supremum over $\nu_A$ on the both sides, which concludes $\tilde V_1^+(t,x,z) \leq V_1^+(t,x,z)$.

By (ii) and (iii), we conclude $V_1^+(t,x,z) = \tilde{V}_1^+(t,x,z)$.

(iv) $V_1^-(t,x) = \tilde V_1^- (t,x) $

Define $\tilde{\mathfrak{B}}_\mu:\mathcal{B}(t)\rightarrow\mathcal{B}(t)$ and its pseudo inverse function $\mathfrak{B}_\mu:\mathcal{B}(t)\rightarrow\mathcal{B}(t)$:
\begin{align}
    & (\tilde{\mathfrak{B}}_\mu(\beta))(s) \coloneqq \begin{cases} \beta\big( \sigma_\mu(s) \big),
    & s\in\text{Range}(\sigma_\mu^{-1}), \\
    \text{any }b \in B, & s\notin\text{Range}(\sigma_\mu^{-1}), \end{cases}
    \label{eq:lemma5_proof4_1}\\
% \end{align}
% \begin{align}
    &(\mathfrak{B}_\mu(\tilde{\beta}))(s) \coloneqq \begin{cases} \tilde{\beta}(\sigma^{-1}_\mu(s)), & s\in[t,\sigma_\mu(T)], \\\text{any }b\in B ,&s\in(\sigma_\mu(T),T],  \end{cases}
    \label{eq:lemma5_proof4_2}
\end{align}
Also, define $\tilde{\mathfrak{C}}_\mu:\Gamma(t)\rightarrow\tilde\Gamma(t)$, where $\tilde\Gamma(t)$ is defined in \eqref{eq:playerA_strategySet_aug}, and its pseudo inverse function $\mathfrak{C}_\mu:\tilde\Gamma(t)\rightarrow\Gamma(t)$:
\begin{align}
    &(\tilde{\mathfrak{C}}_\mu(\gamma)) [\tilde\beta , \mu] (s) =\begin{cases} \gamma [\mathfrak{B}_\mu(\tilde\beta)]\big( \sigma_\mu(s) \big), & s\in\text{Range}(\sigma_\mu^{-1}),\\\text{any } a\in A, & s\notin\text{Range}(\sigma_\mu^{-1}),\end{cases}
    \label{eq:lemma5_proof4_3}\\
    &({\mathfrak{C}_\mu}(\tilde\gamma)) [\beta ] (s) = \begin{cases} \tilde\gamma [\tilde{\mathfrak{B}}_\mu(\beta),\mu]\big( \sigma^{-1}_\mu(s) \big), & s\in[t,\sigma_\mu(T)],\\ \text{any } a\in A, &s\in(\sigma_\mu(T),T].\end{cases}
    \label{eq:lemma5_proof4_4}
\end{align}
% for $\beta\in\mathcal{B}(t)$ and $\tilde{\beta}_{d}\in\mathcal{B}_{d}(t)$.
These definitions satisfy the following properties: for any $\mu\in\mathcal{M}(t)$,
\begin{align}
    % \begin{split}
    %     &(\tilde{\mathfrak{A}}_\mu(\mathfrak{A}_\mu(\tilde\alpha)))(s)= \tilde{\alpha}(s),\\ &(\tilde{\mathfrak{D}}_\mu(\mathfrak{D}_\mu(\tilde\delta)))[\tilde\alpha](s) = \tilde\delta [\tilde\alpha](s), 
    % \end{split}
    % \quad \text{for } s\in\textrm{Range}(\sigma^{-1}_\mu)\\
    % \begin{split}
    %     &(\mathfrak{A}_\mu(\tilde{\mathfrak{A}}_\mu(\alpha)))(s)= {\alpha}(s),\\ &(\mathfrak{D}_\mu(\tilde{\mathfrak{D}}_\mu(\delta)))[\alpha](s) = \delta [\alpha](s), 
    % \end{split}
    % \quad \text{for } s\in[t,\sigma_\mu(T)],\\
    & \big\{ \beta={\mathfrak{B}}_\mu(\tilde\beta)~|~\tilde\beta\in\mathcal{B}(t) \big\} = \mathcal{B}(t),
    \label{eq:lemma5_proof4_5}\\
    & \big\{ \gamma={\mathfrak{C}}_\mu(\tilde\gamma)~|~\tilde\gamma\in\tilde\Gamma(t) \big\} = \Gamma(t).
    \label{eq:lemma5_proof4_6}
\end{align}

Consider $(\tilde{\mathrm{x}},\tilde{\mathrm{z}})$ solving \eqref{eq:dynamics_comb_TI} for $(\tilde\gamma[\tilde\beta,\mu],\tilde\beta,\mu)$,  $(\mathrm{x},\mathrm{z})$ solving \eqref{eq:dynamics_comb} for $(\mathfrak{C}_\mu(\tilde\gamma)[\mathfrak{B}_\mu(\tilde\beta)],\mathfrak{B}_\mu(\tilde\beta))$, and 
$(\mathrm{x}_1,\mathrm{z}_1)$ solving \eqref{eq:dynamics_comb} for $(\gamma[\beta],\beta)$.
\begin{align}
    & \tilde{V}_1^- (t,x,z) = \inf_{\tilde\gamma\in\tilde\Gamma(t)}\sup_{\tilde\beta\in\mathcal{B}(t),\mu\in\mathcal{M}(t)}\max\big\{\max_{s\in[t,T]}c(\tilde{\mathrm{x}}(s)), \notag\\
    &\quad\quad\quad\quad\quad\quad\quad\quad\quad\quad\quad\quad\quad\quad\quad g(\tilde{\mathrm{x}}(T))-\tilde{\mathrm{z}}(T)\big\}\notag\\
    & = \inf_{\tilde\gamma\in\tilde\Gamma(t)}\sup_{\tilde\beta\in\mathcal{B}(t),\mu\in\mathcal{M}(t)}\max\big\{\max_{s\in[t,\sigma_\mu(T)]}c(\mathrm{x}(s)),\notag\\
    &\quad\quad\quad\quad\quad\quad\quad\quad\quad\quad\quad g(\mathrm{x}(\sigma_\mu(T)))-\mathrm{z}(\sigma_\mu(T))\big\}
    \label{eq:lemma5_proof4_7}\\
    & = \inf_{\gamma\in\Gamma(t)}\sup_{\beta\in\mathcal{B}(t),\mu\in\mathcal{M}(t)}\max\big\{\max_{s\in[t,\sigma_\mu(T)]}c(\mathrm{x}_1(s)),\notag\\
    &\quad\quad\quad\quad\quad\quad\quad\quad\quad\quad\quad g(\mathrm{x}_1(\sigma_\mu(T)))-\mathrm{z}_1(\sigma_\mu(T))\big\}.
    \label{eq:lemma5_proof4_8}
\end{align}
\eqref{eq:lemma5_proof4_7} is by \eqref{eq:lemma5_proof0_1} and \eqref{eq:lemma5_proof0_2}, and \eqref{eq:lemma5_proof4_8} is by \eqref{eq:lemma5_proof4_5} and \eqref{eq:lemma5_proof4_6}.
In the term in \eqref{eq:lemma5_proof4_8}, $\mu$ only controls the terminal time ($\sigma_\mu(T)$), hence, the supremum over $\mu$ can be converted to the maximum over $\tau$, which concludes $V_1^-(t,x,z)= \tilde V_1^-(t,x,z)$.

\qed

% %%%%%%%%%%%%%%%%%%%%%%%%%%%%%%%%%%%%%%%%%%%%%%%%%%%%%%%%%%%%%%%%%%%%%%%%%%%%%%%%%%%%%%%%%%%%%%%%%%%%%%%%%%%%%%%%%%%%%%%%%%%%%%%%%%%%%%%%%%%%%%%%%%%%%%%%%%%%%%%%%%%%%
% \subsection{Proof of Lemma \ref{lemma:equiv_V1}}
% \label{appen:lemma_equiv_V1}

% Is there any way to prove this easily?
% Seems not...

%%%%%%%%%%%%%%%%%%%%%%%%%%%%%%%%%%%%%%%%%%%%%%%%%%%%%%%%%%%%%%%%%%%%%%%%%%%%%%%%%%%%%%%%%%%%%%%%%%%%%%%%%%%%%%%%%%%%%%%%%%%%%%%%%%%%%%%%%%%%%%%%%%%%%%%%%%%%%%%%%%%%%
\subsection{Proof of Theorem \ref{thm:HJeq_Prob1_TI}}
\label{appen:thm_HJeq_Prob1_TI}
\noindent
\textbf{Proof.} The terminal value is derived by substituting $T$ for $t$ in \eqref{eq:game_V1_upper_finiteHorizon} or \eqref{eq:game_V1_lower_finiteHorizon}:
\begin{align}
    V_1^\pm (T,x,z) = \max\{c(T,x),g(T,x)-z\}
\end{align}
for all $(x,z)\in\R^n\times\R$.

(i) \cite{altarovici2013general} has presented the HJ equation for state-constrained problems, in which the terminal time is fixed.
By applying the HJ equation in \cite{altarovici2013general} to $V_1^+$, 
\begin{align}
% \begin{split}
    0=\max\big\{ c(x) &- V_1^+ , V^+_{1,t} - \tilde{H}_1^+(x,z,D_x V_1^+, D_z V_1^+) \big\},
    % \\
    %  & - \max_{a\in A} \min_{\substack{b\in B \\ b_d\in [0,1]}} \big[ -D_x V_1^+ \cdot f(x,a,b)b_d + D_z V_1^+L(x,a,b)b_d \big] \bigg\}\notag
% \end{split}
    \label{eq:thm_TI1_proof1}
\end{align}
where 
\begin{align}
    \tilde{H}_1^+(x,z,p,q)\coloneqq \max_{a\in A}\min_{\substack{b\in B \\ b_d\in [0,1]}} -p\cdot f(x,a,b)b_d + qL(x,a,b)b_d.
    \label{eq:thm_TI1_proof2}
\end{align}
% for $(x,z,p,q)\in\R^n\times\R\times\R^n\times\R$.

Since, for all $a\in A , b\in B$, the term $-p\cdot f(x,a,b)b_d + qL(x,a,b)b_d$ is minimized at $b_d=0$ or 1,
\begin{align}
\begin{split}
    \tilde{H}_1^+(x,z,p,q)& = \max_{a\in A}  \min\big\{0,\\
    &\min_{b\in B}-p\cdot f(x,a,b) + qL(x,a,b) \big\}.
\end{split}
\end{align}
Also, $0$ does not depend on $a$, thus, the maximum over $a$ operation can move into the minimum operation:
\begin{align}
\begin{split}
    \tilde{H}_1^+(x,z,p&,q) =   \min\{0,\bar{H}^+(x,z,p,q)\},
\end{split}
\label{eq:thm_TI1_proof3}
\end{align}
where $\bar{H}^+$ is defined in \eqref{eq:def_Hambar_+}.
By applying \eqref{eq:thm_TI1_proof3} to \eqref{eq:thm_TI1_proof1}, \eqref{eq:HJeq1_TI} is proved for $V_1^+$.

(ii) By applying \cite{altarovici2013general} to $V_1^-$, 
\begin{align}
% \begin{split}
    0=\max\big\{ c(x) &- V_1^- , V^-_{1,t} - \tilde{H}_1^-(x,z,D_x V_1^-, D_z V_1^-) \big\},
    \label{eq:thm_TI1_proof4}
\end{align}
where 
\begin{align}
    \tilde{H}_1^-(x,z,p,q)\coloneqq \min_{\substack{b\in B \\ b_d\in [0,1]}}\max_{a\in A} -p\cdot f(x,a,b)b_d + qL(x,a,b)b_d.
    \label{eq:thm_TI1_proof5}
\end{align}
% for $(x,z,p,q)\in\R^n\times\R\times\R^n\times\R$.
Since $b_d\in[0,1]$ is non-negative,
\begin{align}
    \tilde{H}_1^-(x,z,p,q)&= \min_{b_d\in[0,1]} b_d\min_{b\in B}\max_{a\in A} [-p\cdot f(x,a,b) \notag \\
    &\quad\quad\quad\quad\quad\quad\quad\quad\quad +qL(x,a,b)]\notag\\
    % &= \min_{b_d\in[0,1]} b_d \bar{H}^- (x,z,p,q),\notag\\
    & = \min\{0, \bar{H}^- (x,z,p,q) \},
    \label{eq:thm_TI1_proof6}
\end{align}
where $\bar{H}^-$ is defined in \eqref{eq:def_Hambar_-}.
\eqref{eq:thm_TI1_proof4} and \eqref{eq:thm_TI1_proof6} prove \eqref{eq:HJeq1_TI} for $V_1^-$.
\qed
%%%%%%%%%%%%%%%%%%%%%%%%%%%%%%%%%%%%%%%%%%%%%%%%%%%%%%%%%%%%%%%%%%%%%%%%%%%%%%%%%%%%%%%%%%%%%%%%%%%%%%%%%%%%%%%%%%%%%%%%%%%%%%%%%%%%%%%%%%%%%%%%%%%%%%%%%%%%%%%%%%%%%
% \subsection{Proof of Lemma \ref{lemma:equiv_V2}}
% \label{appen:lemma_equiv_V2}

% \textbf{Proof.} hi
%%%%%%%%%%%%%%%%%%%%%%%%%%%%%%%%%%%%%%%%%%%%%%%%%%%%%%%%%%%%%%%%%%%%%%%%%%%%%%%%%%%%%%%%%%%%%%%%%%%%%%%%%%%%%%%%%%%%%%%%%%%%%%%%%%%%%%%%%%%%%%%%%%%%%%%%%%%%%%%%%%%%%
% \subsection{Proof of Theorem \ref{thm:HJeq_Prob2_TI}}
% \label{appen:thm_HJeq_Prob2_TI}

% %%%%%%%%%%%%%%%%%%%%%%%%%%%%%%%%%%%%%%%%%%%%%%%%%%%%%%%%%%%%%%%%%%%%%%%%%%%%%%%%%%%%%%%%%%%%%%%%%%%%%%%%%%%%%%%%%%%%%%%%%%%%%%%%%%%%%%%%%%%%%%%%%%%%%%%%%%%%%%%%%%%%%
% \subsection{Proof of Theorem \ref{thm:Hjeq_prob2_opt_TI}}
% \label{appen:thm_Hjeq_prob2_opt_TI}